\documentclass{article}
\usepackage{amsthm, amssymb, amsmath}
\usepackage{hyperref}

\newtheorem{theorem}{Theorem}[section]

\newtheorem{lemma}[theorem]{Lemma}

\newtheorem{proposition}[theorem]{Proposition}
\newtheorem{corollary}[theorem]{Corollary}

\newcommand{\ip}[2]{\left\langle #1, #2 \right\rangle}

\newcommand{\R}{{\mathbb R}}
\newcommand{\Rp}{{\mathbb R}_+}
\newcommand{\lb}{\left(}
\newcommand{\rb}{\right)}

\renewcommand{\O}[1]{\mathrm{O}\lb #1\rb}
\huge
\include{now_macro}
\begin{document}
\title{The Quenching Problem in the Nonlinear Heat Equations}
\author{Gang Zhou\thanks{Supported by NSERC under Grant NA7901.}}\maketitle
\centerline{Department of Mathematics, University of Toronto, Toronto, Canada, M5S 2E4}
%\fixNumberingInArticle
\setlength{\leftmargin}{.1in}
\setlength{\rightmargin}{.1in} %\centerline{\small{$^{1}$Department
%of Mathematics, University of Toronto, Toronto, Canada}}
%\NowFootNum
%\fixNumberingInArticle
\normalsize \vskip.1in \setcounter{page}{1}
\setlength{\leftmargin}{.1in} \setlength{\rightmargin}{.1in}
\section*{Abstract} In this paper we study the quenching problem in nonlinear heat equations
 with power nonlinearities. For nonlinearities of power $p<0$ and for an open set of slowly varying initial conditions we
prove that the solutions will collapse in a finite time. We find the
collapse profile and estimate the remainder.
\section{Introduction} In this paper we study the problem of collapse of positive solutions for the nonlinear heat equation
\begin{equation}\label{eq:ORIG}
\begin{array}{lll}
\partial_{t}u&=&\partial_{x}^{2}u+[\text{sign}(p-1)]|u|^{p-1}u\\
u(x,0)&=&u_{0}(x).
\end{array}
\end{equation}
For $p>1$ and for suitable initial conditions $u_{0}$ the solutions
of (~\ref{eq:ORIG}) blowup in finite time (see ~\cite{DGSW,GK2, GK3,
HV92, MR1164066, Mer,Vel,MR1230711, MR1317705,MR1433084,BrKu, MZ3,
MZ4, MZ5}). When $p<1$ one expects the solution to collapse in
finite time for certain initial conditions. It is the second case
which is of the interest in this paper. The $p<1$ problem arises in
the study of the quenching problem in combustion theory ~\cite{Guo1,
Guo2} and vortex reconnection ~\cite{MZ2}. Moreover it presents a
simple mathematical model for the neck pinching problem of mean
curvature flow and Ricci flow ~\cite{MR1737226, SK, SS, Huis1,
Huis2, EvSp} and for collapse in the Keller-Segal problem of
chemotaxis ~\cite{MR1709861}.

We say that a positive solution $u$ collapses at time $t^*$ if
$u(\cdot,t)\geq c(t)>0$ for $t<t^*$ and some positive scalar
function $c(t)$ and $u\rightarrow 0$ as $t\rightarrow t^*$ on some
set $x\in \mathcal{S}\subset \mathbb{R}$. Because of its application
in combustion theory the problem of collapse of solutions is
referred to as the quenching problem. In this paper, for technical
reason, we limit ourselves to the case $p<0$ and $u_{0}\geq c_{0}>0$
for some constant $c_{0}$.

The problem of quenching for (~\ref{eq:ORIG}) with $p=-1$ on bounded
domain was studied first in ~\cite{KH}, where a sufficient condition
for collapsing is found. Later Huisken ~\cite{Huis2} proved that if
$\partial_{x}^{2}u_{0}(x)-u^{p}_{0}(x)\leq 0$, then the solution
collapses in finite time. Merle and Zaag proved in ~\cite{MZ2} that
there exists initial condition $u_{0}(x)$ such that the solution
$u(x,t)$ collapses in finite time $t^*$ and
\begin{equation}\label{eq:convergence}
\lim_{t\rightarrow t^*}\|(t^*-t)^{-\frac{1}{1-p}}u(x((t^*-t)|ln(t^*-t)|)^{1/2},t)-
(1-p-\frac{(1-p)^{2}}{4p}x^{2})^{\frac{1}{1-p}}\|_{\infty}=0.
\end{equation}
For the neck pinching problem in mean curvature flow, Huisken proved
in ~\cite{Huis2} a result weaker than (~\ref{eq:convergence}) holds
on a bounded space domain. For other related works we refer to
~\cite{Guo1, Guo2, SS, Huis1}. The starting point in these works is
to study the rescaled function
$(t^*-t)^{-\frac{1}{2}}u(x(t^*-t)^{\frac{1}{2}},t)$ using the
technique of Sturm Liouville theorem for linear parabolic equations,
as used in ~\cite{GK1,GK2}.

The origin of scaling and asymptotic in (~\ref{eq:convergence}) lies
in the following key properties of Equation (~\ref{eq:ORIG}):
\begin{enumerate}
 \item (~\ref{eq:ORIG}) is invariant
 with respect to the scaling transformation,
 \begin{equation}\label{rescale}
 u(x,t)\rightarrow \lambda^{\frac{2}{1-p}}
u(\lambda^{-1} x,\lambda^{-2} t)
\end{equation} for any constant $\lambda>0,$ i.e. if
$u(x,t)$ is a solution, so is $\lambda^{\frac{2}{1-p}}u(\lambda^{-1}
x,\lambda^{-2}t).$
 \item (~\ref{eq:ORIG}) has $x-$independent (homogeneous)
 solutions:
\begin{equation}\label{eq:homegeous}
 u_{hom}=[u_{0}^{-p+1}-|p-1|t]^{-\frac{1}{p-1}}.
\end{equation} These solutions collapse (blow up) in finite time
 $t^* = ( |p-1|u_0^{p-1})^{-1}$.
\end{enumerate}

In what follows we use the notation $f\lesssim g$ for two functions
$f$ and $g$ satisfying $f\leq Cg$ for some universal constant $C$.
We will also deal, without specifying it, with weak solutions of
Equation (~\ref{eq:ORIG}) in some appropriate sense. These solutions can
be shown to be classical for $t>0.$

In our paper we consider Equation (~\ref{eq:ORIG}) in $~\mathbb{R}$
with $p<0$ and with the initial conditions $u_{0}$ even, bounded
below by a positive constant and having a local minimum at $0$
modulo small fluctuation. We prove that there exists a time
$t^{*}<\infty$ such that $u(x,t)$ collapses at time $t^*$; moreover
there exist $C^{1}$ functions $\lambda(t)$, $b(t)$, $c(t)$ and
$\eta(x,t)$ such that
\begin{equation}\label{eq:DecomU}
u(x,t)=\lambda^{\frac{2}{1-p}}(t)[(\frac{1-p+b(t)\lambda^{-2}(t)x^{2}}{2c(t)})^{\frac{1}{1-p}}+\eta(x,t)]\end{equation}
with
\begin{equation}\label{eq:Remain}
\|\langle \lambda^{-1}(t)x\rangle^{-3}\eta(x,t)\|_{\infty}\lesssim
b^{3/2}(t).
\end{equation} Furthermore the scalar functions $\lambda(t), \ b(t)$ and $c(t)$
satisfy the estimates
\begin{equation}\label{eq:para}
\begin{array}{lll}
\lambda(t)&=&\lambda(0)(t^*-t)^{\frac{1}{2}}(1+o(1)),\ \lambda(0)=(2c_{0}+\frac{2}{1-p}b_{0})^{-1/2};\\
& &\\
b(t)&=&\frac{(p-1)^{2}}{4pln|t^*-t|}(1+O(\frac{1}{|ln|t^*-t|}|^{1/2}));\\
& &\\
c(t)&=&\frac{1}{2}+\frac{1-p}{4p|ln|t^*-t||}(1+O(\frac{1}{ln|t^*-t|}))
\end{array}
\end{equation} where $o(1)\rightarrow 0$ as $t\rightarrow t^*.$

Before stating the main theorem we define a function $g(y,b_{0})$ as
$$g(y,b_{0}):=\left\{
\begin{array}{lll}
(\frac{1-p}{2})^{\frac{1}{1-p}}\ \ \ \ \ \ \text{if}\ b_0 y^{2}\leq 4(1-p)\\
(2(1-p))^{\frac{1}{1-p}}\ \text{if}\ b_0 y^{2}>4(1-p)
\end{array}
\right.$$ and define the constant $q$ which will be used throughout
the rest of paper $$q:=\min\{\frac{4}{1-p},\
\frac{2(2-p)}{(1-p)^{2}},\ 1\}.$$ The following is the main result
of our paper.
\begin{theorem}\label{maintheorem} Assume the initial datum
$u_{0}(x)$ in (~\ref{eq:ORIG}) is even and satisfies the estimates
\begin{equation}\label{eq:INI2}
\begin{array}{lll}
\|\langle x\rangle^{-n}[u_{0}(x)-(\frac{1-p+
b_{0}x^{2}}{2c_{0}})^{\frac{1}{1-p}}]\|_{\infty}&\leq&
\delta_{0}b_{0}^{\frac{n}{2}},
\end{array}
\end{equation} $n=2,3,\ u_{0}(x)\in \langle x\rangle^{\frac{2}{1-p}} L^{\infty}$
and $u_{0}(x)\geq (2c_{0}+\frac{2
b_{0}}{1-p})^{\frac{1}{p-1}}g((2c_{0}+\frac{2
b_{0}}{1-p})^{1/2}x,b_{0})$ for some $1/2\leq c_{0}\leq 2$, and (~\ref{eq:INI2}) with $n=q$ if
$p<-1$. Then there exists a constant $\delta$
such that if $\delta_{0},\ b_{0}\leq \delta$
 then there exists a finite
time $0<t^{*}< \infty$ such that
$\|\frac{1}{u(\cdot,t)}\|_{\infty}<\infty$ if $t<t^*$ and
$$\|\frac{1}{u(\cdot,t)}\|_{\infty}\rightarrow \infty\ \text{as}\
t\rightarrow t^*.$$ Moreover there exist $C^1$ functions
$\lambda(t)$, $b(t)$, $c(t)$ and $\eta(x,t)$ such that $u(x,t)$
satisfies the estimates in (~\ref{eq:DecomU})- (~\ref{eq:para}).
\end{theorem}
The proof of this theorem is given in Section ~\ref{SecMain}. This
theorem shows the collapse at $0$ for a certain neighborhood of the
homogeneous solution, \eqref{eq:homegeous}, and it provides a
detailed description of the leading term and an estimate of the
remainder in $\langle x\rangle^{3}L^\infty$. In fact, we have not
only the asymptotic expressions for the parameters $b$ and $c$
determining the leading term and the size of the remainder, but also
dynamical equations for these parameters:
\begin{align}
 b_\tau&=\frac{4 p}{(p-1)^{2}}
b^2+c^{-1}c_{\tau}
b+{\mathcal R}_b(\eta,b,c),\\
c^{-1}c_\tau&=2(\frac{1}{2}-c)-\frac{2}{p-1}b+{\mathcal
R}_c(\eta,b,c),
\end{align}
where $\tau$ is a `collapse' time related to the original time $t$
as $\tau(t):=\int_{0}^{t}\lambda^{-2}(s)ds$ and the remainders have
the estimates
\begin{equation}
{\mathcal R}_b(\eta,b,c),{\mathcal R}_c(\eta,b,c) =\O{
b^3+[|c-\frac{1}{2}|+|c_\tau|]b^2+|b_\tau|
b+b\|\eta(\cdot,t)\|_{X}+\|\eta(\cdot,t)\|_{X}^{2}+\|\eta(\cdot,t)\|_{X}^{2-p}}
\end{equation}
with the norm $\|\eta(\cdot,t)\|_{X}:=\|\langle
\lambda^{-1}(t)x\rangle^{-3}\eta\|_{\infty}.$

This paper is organized as follows.  In Section ~\ref{SEC:LocalWell}
we prove the local well-posedness of Equation (~\ref{eq:ORIG}). In
Sections \ref{SEC:BlowVaria}-\ref{SEC:RePara} we present some
preliminary derivations and motivations for our analysis. In Section
\ref{SEC:ApriEst}, we formulate a priori bounds on solutions to
\eqref{eq:ORIG} which are proven in Sections \ref{SEC:EstB}-
~\ref{SEC:EstM3}. We use these bounds and a lower bound proved in
Section ~\ref{SEC:LowerBound} in Section \ref{SecMain} to prove our
main result, Theorem \ref{maintheorem}. In Section
\ref{SEC:ParaMetrization} we lay the ground work for the proof of
the a priori bounds of Section \ref{SEC:ApriEst}, in particular by
using a Lyapunov-Schmidt-type argument we derive equations for the
parameters $a$, $b$ and $c$ and fluctuation $\eta$.

\section*{Acknowledgement}
The author would like to thank his Ph.D advisor, Professor I.M.
Sigal, for many discussion and important suggestions.
\section{Local Well-posedness of
(~\ref{eq:MCF})}\label{SEC:LocalWell}
% Proofread\\
In this section we prove the local well-posedness of
(~\ref{eq:ORIG}) for $p<0$ in spaces used in this paper. Since in
what follows we are dealing with $p<0$ and $u_{0}(x)>0$ we restate
(~\ref{eq:ORIG}) as
\begin{equation}\label{eq:MCF}
\begin{array}{lll}
\partial_{t}u&=&\partial^{2}_{x}u-u^{p},\ p<0\\
u(x,0)&=&u_{0}(x)>0.
\end{array}
\end{equation}

\begin{theorem}\label{THM:WellPose}
If $u_{0}(x)\in \langle x\rangle^{\frac{2}{1-p}} L^{\infty}$ and
$u_{0}(x)\geq \kappa_{0}>0,$ then there exist a function $f(x,t)\in
\langle x\rangle^{\frac{2}{1-p}} L^{\infty}$ for any time $t<\infty$
and a time $\delta(\kappa_{0},\|\langle
x\rangle^{-\frac{2}{1-p}}u_{0}\|_{\infty})$ such that for any time
$t_{0}\leq t\leq t_{0}+\delta$, (~\ref{eq:MCF}) has a unique
solution $u(\cdot,t)\in \langle x\rangle^{\frac{2}{1-p}}L^{\infty}$
with $u(x,0)=u_{0}(x)$ and $f(x,t)\geq u(x,t)\geq
\frac{1}{2}\kappa_{0}.$

Moreover, if $t_*$ is the supremum of such
$\delta(\kappa_{0},\|\langle
x\rangle^{-\frac{2}{1-p}}u_{0}\|_{\infty})$ then either
$t_{*}=\infty$ or $\|\frac{1}{u(\cdot,t)}\|_{\infty}\rightarrow
\infty$ as $t\rightarrow t_*.$
\end{theorem}
\begin{proof}
We transform (~\ref{eq:MCF}) as
\begin{equation}\label{eq:subsolu}
u(x,t)=e^{t\partial_{x}^{2}}u_{0}(x)-\int_{0}^{t}e^{(t-s)\partial_{x}^{2}}u^{p}(x,s)ds.
\end{equation}
Define a new function $u_{1}(x,t)$ by
$$
u_{1}(x,t):=u(x,t)-f(x,t) $$ with
$f(x,t):=e^{t\partial_{x}^{2}}u_{0}(x).$ Then (~\ref{eq:subsolu})
becomes
\begin{equation}\label{eq:u1}
u_{1}(x,t)=-\int_{0}^{t}e^{(t-s)\partial_{x}^{2}}[f(x,s)+u_{1}(x,s)]^{p}ds.
\end{equation}
In the next we use the fix point argument to prove the existence and
the uniqueness of the solution $u_{1}$ to (~\ref{eq:u1}), hence
those of $u$ to complete the proof.

We start with proving $f(x,t)\geq \kappa_{0}>0$ and $\langle
x\rangle^{-\frac{2}{1-p}}f(\cdot,t)\in L^{\infty}$. It is well known
that the integral kernel of $e^{t\partial_{x}^{2}}$ is
$\frac{1}{\sqrt{\pi t}}e^{-\frac{(x-y)^2}{t}}>0$, consequently
\begin{equation}\label{eq:kernel}
f(x,t)\geq e^{t\partial_{x}^{2}}\kappa_{0}= \kappa_{0}.
\end{equation} Moreover we claim that $\langle x\rangle^{-\kappa}e^{t\partial_{x}^{2}}\langle x\rangle^{\kappa}<\infty$
for any $2\geq \kappa\geq 0$. Hence by the fact $2>\frac{2}{1-p}>0$ we
have $f(\cdot,t)\in \langle x\rangle^{\frac{2}{1-p}}L^{\infty}.$ In
the following we prove the claim for $\kappa=0$ and $\kappa=2$, the general
case follows by interpolation between them. It is
easy to prove $\kappa=0$ by the fact $e^{t\partial_{x}^{2}}1=1$; for
$\kappa=2$ the
fact $y^{2}\leq 2(x-y)^{2}+2x^{2}$ yields
$$\langle
x\rangle^{-2}e^{t\partial_{x}^{2}}x^{2}=\langle
x\rangle^{-2}\int\frac{1}{\sqrt{\pi
t}}e^{-\frac{(x-y)^{2}}{t}}y^{2}dy\leq 2\langle x\rangle^{-2}\int
\frac{1}{\sqrt{\pi
t}}e^{-\frac{(x-y)^{2}}{t}}[(x-y)^{2}+x^{2}]dy=\langle
x\rangle^{-2}[2x^{2}+td_{0}]<\infty$$ with the constant $d_{0}:=\frac{2}{\sqrt{\pi }}\int
e^{-x^{2}}x^{2}d x.$ Thus we finish proving
the properties of $f(\cdot,t).$

By the integral kernel we have $f(x,t)\rightarrow u_{0}(x)$ as
$t\rightarrow 0^+.$ Moreover by using the contraction lemma on
(~\ref{eq:u1}) it is not hard to prove that there exists a time $\delta(\kappa_{0},\|\langle x\rangle^{-\frac{2}{1-p}}u_{0}\|_{\infty})$ such that for time $t\in [0,\delta]$ there exists a unique
bounded negative solution $u_{1}$ such that
\begin{equation}\label{eq:ESTu1}
\|u_{1}(\cdot,t)\|_{\infty}\leq \frac{\kappa_{0}}{2}.
\end{equation}

Thus $u(x,t)=f(x,t)+u_{1}(x,t)\geq \frac{\kappa_{0}}{2}$ has all the properties in
the proposition, thus the proof is complete.
\end{proof}

\section{Blow-Up Variables and Almost
Solutions}\label{SEC:BlowVaria}
%%%%%%%%%%%%%%%%%%%%%%%%%%%%%%%%%%%%%%%%%%%%%%%%%%%%%%%%%%%%%%%%%%%%%%%%
%%%%%%%%%%%%%%%%%%%%%%%%%%%%%%%%%%%%%%%%%%%%%%%%%%%%%%%%%%%%%%%%%%%%%%%%
% Proofread\\
In this section we pass from the original variables $x$ and $t$ to
the blowup variables $y:=\lambda^{-1}(t) (x-x_{0}(t))$ and
$\tau:=\int_{0}^{t} \lambda^{-2}(s) ds$.  The point here is that we
do not fix $\lambda(t)$ and $x_0(t)$ but consider them as free
parameters to be found from the evolution of \eqref{eq:MCF}. Assume
for simplicity that $u_0$ is even with respect to $x=0$. In this
case $x_{0}$ can be taken to be $0$. Suppose $u(x,t)$ is a solution
to \eqref{eq:MCF} with an initial condition $u_{0}(x)$. We define
the new function
\begin{equation}\label{eq:definev}
v(y,\tau):=\lambda^{\frac{2}{p-1}}(t)u(x,t)
\end{equation} with $y:=\lambda^{-1}(t) x$ and $\tau:=\int_{0}^{t}\lambda^{-2}(s)ds.$
The function $v$ satisfies the equation
\begin{equation}\label{eqn:BVNLH}
v_{\tau}=(\partial_y^2-a y\partial_y+\frac{2 a}{1-p}) v-v^{p}.
\end{equation}
where $a:=-\lambda \partial_{t}\lambda$.  The initial condition is
$v(y,0)=\lambda^{\frac{2}{p-1}}_{0} u_{0}(\lambda_{0}y)$, where
$\lambda_0$ is the initial condition for the scaling parameter
$\lambda$.

%%%%%%%%%%%%%%%%%%%%%%%%%%%%%%%%%%%%%%%%%%%%%%%%%%%%%%%%%%%%%%%%%%%%%%%%
%%%%%%%%%%%%%%%%%%%%%%%%%%%%%%%%%%%%%%%%%%%%%%%%%%%%%%%%%%%%%%%%%%%%%%%%
%\section{Almost Solutions to the Nonlinear Heat Equation in Blowup
%Variables} \label{sec:AlmostSolutions}
%%%%%%%%%%%%%%%%%%%%%%%%%%%%%%%%%%%%%%%%%%%%%%%%%%%%%%%%%%%%%%%%%%%%%%%%
%%%%%%%%%%%%%%%%%%%%%%%%%%%%%%%%%%%%%%%%%%%%%%%%%%%%%%%%%%%%%%%%%%%%%%%%
If the parameter $a$ is a constant, then \eqref{eqn:BVNLH} has the
following homogeneous, static (i.e. $y$ and $\tau$-independent)
solutions
\begin{equation}
v_a:=( \frac{1-p}{2 a})^\frac{1}{1-p}.
\end{equation}
In the original variables $t$ and $x$, this family of solutions
corresponds to the homogeneous solution \eqref{eq:homegeous} of the
nonlinear heat equation with the parabolic scaling $\lambda^{2}=2
a(T-t)$, where the collapsing time,
$T:=\left[u_0^{p-1}(1-p)\right]^{-1}$, is dependent on $u_0$, the
initial value of the homogeneous solution $u_{hom}(t)$.

If the parameter $a$ is $\tau$ dependent but $|a_\tau|$ is small,
then the above solutions are good approximations to the exact
solutions. Another approximation is the solution of $ a y
v_y+\frac{2 a}{p-1} v+v^p=0$, obtained from \eqref{eqn:BVNLH} by
neglecting the $\tau$ derivative and second order derivative in $y$.
This equation has the general solution
\begin{equation}
v_{a b}:=(\frac{1-p + b y^2}{2 a})^\frac{1}{1-p}
\end{equation}
for all $b\in\mathbb{R}$. In what follows we take $b\ge 0$ so that
$v_{a b}$ is nonsingular.  Note that $v_{a,0}=v_a$.

%%%%%%%%%%%%%%%%%%%%%%%%%%%%%%%%%%%%%%%%%%%%%%%%%%%%%%%%%%%%%%%%%%%%%%%%
%%%%%%%%%%%%%%%%%%%%%%%%%%%%%%%%%%%%%%%%%%%%%%%%%%%%%%%%%%%%%%%%%%%%%%%%
%\appendix
%%%%%%%%%%%%%%%%%%%%%%%%%%%%%%%%%%%%%%%%%%%%%%%%%%%%%%%%%%%%%%%%%%%%%%%%
%%%%%%%%%%%%%%%%%%%%%%%%%%%%%%%%%%%%%%%%%%%%%%%%%%%%%%%%%%%%%%%%%%%%%%%%

%%%%%%%%%%%%%%%%%%%%%%%%%%%%%%%%%%%%%%%%%%%%%%%%%%%%%%%%%%%%%%%%%%%%%%%%
%%%%%%%%%%%%%%%%%%%%%%%%%%%%%%%%%%%%%%%%%%%%%%%%%%%%%%%%%%%%%%%%%%%%%%%%
\section{``Gauge" Transform}
\label{sec:Gauge}
%%%%%%%%%%%%%%%%%%%%%%%%%%%%%%%%%%%%%%%%%%%%%%%%%%%%%%%%%%%%%%%%%%%%%%%%
%%%%%%%%%%%%%%%%%%%%%%%%%%%%%%%%%%%%%%%%%%%%%%%%%%%%%%%%%%%%%%%%%%%%%%%%
% Proofread\\
We assume that the parameter $a$ depends slowly on $\tau$ and treat
$|a_\tau|$ as a small parameter in a perturbation theory for
Equation \eqref{eqn:BVNLH}.
 In order to convert the global non-self-adjoint operator $a y\partial_y$
appearing in this equation into a more tractable local and
self-adjoint operator we perform a gauge transform. Let
\begin{equation}\label{eq:definew}
w(y,\tau):=e^{-\frac{a y^2}{4}}v(y,\tau).
\end{equation}
Then $w$ satisfies the equation
\begin{equation} \label{eqn:w}
w_\tau=( \partial_y^2-\frac{1}{4}\omega^2
y^2-(\frac{2}{p-1}-\frac{1}{2}) a ) w-e^{\frac{a}{4}(p-1)y^2}w^p,
\end{equation}
where $\omega^2=a^2+a_\tau$.  The approximate solution $v_{a b}$ to
\eqref{eqn:BVNLH} transforms to $v_{a b c}$ where $v_{a b c}:=v_{
c,b}e^{-\frac{a y^2}{4}}$. Explicitly
\begin{equation}
v_{a b c}:=(\frac{1-p+b y^2}{2 c})^{\frac{1}{1-p}}e^{-\frac{a
y^2}{4}}.
\end{equation}
%%%%%%%%%%%%%%%%%%%%%%%%%%%%%%%%%%%%%%%%%%%%%%%%%%%%%%%%%%%%%%%%%%%%%%%%
%%%%%%%%%%%%%%%%%%%%%%%%%%%%%%%%%%%%%%%%%%%%%%%%%%%%%%%%%%%%%%%%%%%%%%%%
\section{Re-parametrization of Solutions}\label{SEC:RePara}
%%%%%%%%%%%%%%%%%%%%%%%%%%%%%%%%%%%%%%%%%%%%%%%%%%%%%%%%%%%%%%%%%%%%%%%%
%%%%%%%%%%%%%%%%%%%%%%%%%%%%%%%%%%%%%%%%%%%%%%%%%%%%%%%%%%%%%%%%%%%%%%%%
In this section we split solutions to \eqref{eqn:w} into the leading
term - the almost solution $v_{a b c}$ - and a fluctuation $\xi$
around it.  More precisely, we would like to parametrize a solution
by a point on the manifold $M_{as}:=\{v_{a b c}\, |\, a,b,c\in \Rp,
b \leq \epsilon,\, a=a(b,c)\}$ of almost solutions and the
fluctuation orthogonal to this manifold (large slow moving and small
fast moving parts of the solution).  Here $a=a(b,c)$ is a twice
differentiable function of $b$ and $c$.  For technical reasons, it
is more convenient to require the fluctuation to be almost
orthogonal to the manifold $M_{as}$. More precisely, we require
$\xi$ to be orthogonal to the vectors $\phi_{0 a}:=(\frac{\alpha}{2\pi})^{\frac{1}{4}}e^{-\frac{a}{4}
y^2}$ and $\phi_{2 a}:=(\frac{\alpha}{8\pi})^{\frac{1}{4}}(1-a y^2)e^{-\frac{a}{4} y^2}$ which are
almost tangent vectors to the above manifold, provided $b$ is
sufficiently small.  Note that $\xi$ is already orthogonal to
$\phi_{1 a}:=(\frac{\alpha}{2\pi})^{\frac{1}{4}}\sqrt{a} y e^{-\frac{a}{4} y^2}$ since our initial
conditions, and therefore, the solutions are even in $x$.

In what follows we fix the relation between $a$ and $c$ as
$$2c(\tau)=a(\tau)+\frac{1}{2}.$$

Define a new function
$V_{a,b}:=(\frac{1-p+by^{2}}{a+\frac{1}{2}})^{\frac{1}{1-p}}$ and a
neighborhood $U_{\epsilon_{0}}$:
\begin{equation*}
U_{\epsilon_0}:=\{v\in L^\infty(\R)\ |\ \|e^{-\frac{1}{8}
y^2}(v-V_{ab})\|_{{\infty}}=o(b)\ \mbox{for some}\ a\in[1/4,1],\
b\in[0,\epsilon_0]\ \}.
\end{equation*}

\begin{proposition}\label{Prop:Splitting}
There exist an $\epsilon_{0}>0$ and a unique $C^1$ functional
$g:U_{\epsilon_0}\rightarrow \mathbb{R}^{+}\times \mathbb{R}^{+}$,
such that any function $v\in U_{\epsilon_0}$ can be uniquely written
in the form
\begin{equation} \label{eqn:split}
v =V_{g(v)} + \eta,
\end{equation}
with $\eta\perp e^{-\frac{a}{4} y^2}\phi_{0 a},\ e^{-\frac{a}{4}
y^2}\phi_{2 a}$ in $L^2(\R)$, $(a,b)=g(v)$. Moreover if
$\|e^{-\frac{1}{8} y^2}(v-V_{a_{0}b_{0}})\|_{{\infty}}=o(b_{0})$
for some $a_{0}$ and $b_{0}$ then
\begin{equation}\label{eq:approx}
|g(v)-(a_{0}, b_{0})|\lesssim
\|e^{-\frac{1}{8}y^{2}}(v-V_{a_{0},b_{0}})\|_{\infty}.
\end{equation}
\end{proposition}
\begin{proof}
 The orthogonality conditions on the fluctuation can be written as
$G(\mu,v)=0$, where $\mu=(a,b)$ and $G:\mathbb{R}^{+}\times
\mathbb{R}^{+}\times
 L^\infty( \R)\rightarrow \R^2$ is defined
as
\begin{equation*}
G(\mu, v):=\left( \begin{array}{c} \ip{V_{\mu}-v}{(\frac{\alpha}{2\pi})^{-\frac{1}{4}}e^{-\frac{ay^{2}}{4}}\phi_{0a}}\\
\ip{V_{\mu}-v}{(\frac{\alpha}{8\pi})^{-\frac{1}{4}}e^{-\frac{ay^{2}}{4}}\phi_{2a}}
\end{array} \right).
\end{equation*}
Here and in what follows, all inner products are $L^2$ inner
products.  Using the implicit function theorem we will prove that
for any $\mu_0:=(a_0, b_0)\in [\frac{1}{4},1]\times (0,\epsilon_{0})
$ there exists a unique $C^1$ function $g:L^\infty\rightarrow
\mathbb{R}^{+}\times \mathbb{R}^{+}$ defined in a neighborhood
$U_{\epsilon_0}\subset L^\infty$ of $V_{\mu_0}$ such that
$G(g(v),v)=0$ for all $v\in U_{\epsilon_0}$.

Note first that the mapping $G$ is $C^1$ and $G(\mu_0, V_{\mu_0})=0$
for all $\mu_0$.  We claim that the linear map $\partial_\mu
G(\mu_0, V_{\mu_0})$ is invertible.  Indeed, we compute
\begin{equation}\label{linearization}
\partial_\mu G(\mu, v)|_{\mu=\mu_{0}}=A_{1}(\mu)+A_{2}(\mu,v)|_{\mu=\mu_{0}}
\end{equation}
where
$$
A_{1}:=\left(
\begin{array}{cc}
\ip{\partial_a V_{
\mu}}{ e^{-\frac{a}{2} y^2}} & \ip{\partial_b V_{\mu}}{ e^{-\frac{a}{2} y^2}}  \\
\ip{\partial_a V_{\mu}}{(1-a y^2)e^{-\frac{a}{2} y^2}}&
\ip{\partial_b V_{\mu}}{(1-a y^2)e^{-\frac{a}{2} y^2}}
\end{array}
\right)
$$ and
\begin{equation*}
A_{2}:= -\frac{1}{4}\left(\begin{array}{cc}
\langle V_{\mu}-v, y^{2}e^{-\frac{a}{2}y^{2}}\rangle& 0 \\
\ip{V_{\mu}-v}{( 5-a y^2) y^2 e^{-\frac{a}{2} y^2}} & 0
\end{array}\right).
\end{equation*}
By the condition in the proposition, we have for $A_{2}$ that
$$\|A_{2}(\mu_{0},v)\|\lesssim |b-b_{0}|+|a-a_{0}|+|b_{0}|.$$  For $b>0$ and
small, we expand the matrix $A_{1}$ in $b$ to get
$A_{1}=G_{1}G_{2}+o(b)$, where the matrices $G_{1}$ and $G_{2}$ are
defined as
$$G_{1}:=\left(
\begin{array}{ccc}
\langle -y^{2}e^{-\frac{ay^{2}}{4}},e^{-\frac{ay^{2}}{4}}\rangle &
\frac{1}{a+\frac{1}{2}}\langle e^{-\frac{ay^{2}}{4}},e^{-\frac{ay^{2}}{4}}\rangle\\
\langle -y^{2}e^{-\frac{ay^{2}}{4}},
(1-ay^{2})e^{-\frac{ay^{2}}{4}}\rangle & 0
\end{array}
\right)$$ and
$$G_{2}:=(\frac{1-p}{a+1/2})^{\frac{1}{1-p}}\frac{1}{p-1}\left(
\begin{array}{lll}
\frac{p-1}{4} & 1\\
1 & 0
\end{array}
\right).$$ Obviously the matrices $G_{1}$ and $G_{2}$ have uniformly
(if $a\in [\frac{1}{4},1]$) bounded inverses.  We claim that this
observation implies our proposition. Indeed, expand $G(\mu,v)$ as
$$G(\mu,v)=G(\mu_{0},v)+\partial_{\mu}G(\mu_{0},v)\left(
\begin{array}{lll}
a-a_{0}\\
b-b_{0}
\end{array}
\right)+O(|b-b_{0}|^{2}+|a-a_{0}|^{2}),$$ provided that $b,
b_{0}>0.$ By the assumption in the proposition we have that
$$
G(\mu_{0},v)=o(b_{0})\ \text{if}\ a_{0}\in[\frac{1}{4},1].
$$ Moreover as we have seen above, the matrix
$\partial_{\mu}G(\mu_{0},v)$ for $v\in U_{\epsilon_{0}}$ has
uniformly bounded inverse. Hence we have that $G(\mu,v)=0$ has a
unique solution $g(v)$ satisfying
$$|g(v)-(a_{0},b_{0})|\lesssim |G(\mu_{0},v)|\lesssim \|e^{-\frac{1}{8}y^{2}}(v-V_{a_{0},b_{0}})\|_{\infty}$$ which is
(~\ref{eq:approx}).
\end{proof}
Recall that $q=\min\{\frac{4}{1-p},\ \frac{2(2-p)}{(1-p)^{2}},\
1\}.$
\begin{proposition}
\label{Prop:SplittingIC} In the notation of Proposition
\ref{Prop:Splitting}, if $\|\langle
y\rangle^{-n}(v-V_{\mu_{0}})\|_{{\infty}}\le
\delta_{0}b_{0}^{\frac{n}{2}}$ with $n=2,3,$ and $n=q$ if $p<-1$,
with $\delta_{0}, b_{0}$ small, then
\begin{equation}\label{eq:vv0}
|g(v)-\mu_{0}|\lesssim \delta_{0}b_{0}^{3/2};
\end{equation}
\begin{equation}
\|\langle y\rangle^{-3}(v-V_{g(v)})\|_{{\infty}}\lesssim \|\langle
y\rangle^{-3}(v-V_{\mu_{0}})\|_{{\infty}}; \label{Ineq:IC}
\end{equation}
\begin{equation}\label{eq:vwithoutweight}
\|\langle y\rangle^{-2}(v-V_{g(v)})\|_{{\infty}}\lesssim
\delta_{0}b_{0}+\delta_{0}b_{0}^{3/2};
\end{equation} and if $p<-1$ then
\begin{equation}\label{eq:vh}
\|\langle y\rangle^{-q}(v-V_{g(v)})\|_{{\infty}}\lesssim
\delta_{0}b_{0}^{q/2}+\delta_{0}b_{0}^{\frac{1+q}{2}}.
\end{equation}
\end{proposition}
\begin{proof}
Equation (~\ref{eq:approx}) implies that
\begin{equation}\label{eq:vmu}
|g(v)-\mu_{0}| \lesssim\|\langle
y\rangle^{-3}(v-V_{\mu_{0}})\|_{{\infty}}
\end{equation} with
$\mu_{0}:=(a_{0},b_{0}).$ This together with the fact $\|\langle
y\rangle^{-3}(v-V_{\mu_{0}})\|_{\infty}\leq \delta_{0}b_{0}^{3/2}$
yields (~\ref{eq:vv0}). Moreover
$$
\begin{array}{lll}
\|\langle y\rangle^{-3}(v-V_{g(v)})\|_{{\infty}}&\leq &\|\langle
y\rangle^{-3}(v-V_{\mu_{0}})\|_{{\infty}}+\|\langle
y\rangle^{-3}(V_{g(v)}-V_{\mu_{0}})\|_{{\infty}}\\
&\lesssim& \|\langle
y\rangle^{-3}(v-V_{\mu_{0}})\|_{{\infty}}+|\mu_{0}-g(v)|\\
&\lesssim & \|\langle y\rangle^{-3}(v-V_{\mu_{0}})\|_{{\infty}}
\end{array}
$$ which is (~\ref{Ineq:IC}).

To prove Equation (~\ref{eq:vwithoutweight}),
%is implied by (~\ref{eq:approx}).
we write
$$\|\langle y\rangle^{-2}(v-V_{g(v)})\|_{{\infty}}\leq \|\langle y\rangle^{-2}(v-V_{\mu_{0}})\|_{{\infty}}+\|\langle y\rangle^{-2}(V_{g(v)}-V_{\mu_{0}})\|_{{\infty}}.$$
By (~\ref{eq:vmu}), we have $\|\langle
y\rangle^{-2}(V_{g(v)}-V_{\mu_{0}})\|_{\infty}\lesssim|g(v)-\mu_{0}|\lesssim
\delta_{0}b_{0}^{3/2}$, hence
$$\|\langle y\rangle^{-2}(V_{g(v)}-V_{\mu_{0}})\|_{{\infty}}\lesssim
|g(v)-\mu_{0}|\lesssim \delta_{0}b_{0}^{3/2}.$$ This together with the fact
$\|\langle y\rangle^{-2}(v-V_{\mu_{0}})\|_{{\infty}}\leq \delta_{0}b_{0}$
completes the proof of (~\ref{eq:vwithoutweight}).

Finally for Equation (~\ref{eq:vh}) we have
$$\|\langle y\rangle^{-q}(v-V_{g(v)})\|_{{\infty}}\leq \|\langle y\rangle^{-q}(v-V_{\mu_{0}})\|_{{\infty}}+\|\langle y\rangle^{-q}(V_{g(v)}-V_{\mu_{0}})\|_{{\infty}}.$$
By its definition we have $1\geq q>\frac{2}{1-p}$ which together
with the definition of $V_{a,b}$ yields
$$
\begin{array}{lll}
\|\langle y\rangle^{-q}(V_{g(v)}-V_{\mu_{0}})\|_{{\infty}}&\lesssim
& |a-a_{0}|+\|\frac{|b-b_{0}||y|^{2-q}}{(1-p+b
y^{2})^{1-\frac{1}{1-p}}}\|_{\infty}\\
&\lesssim &|a-a_{0}|+|b-b_{0}|b_{0}^{\frac{q}{2}-1}.
\end{array}
$$ By the estimate of $|a-a_{0}|+|b-b_{0}|$ above
we have $$\|\langle
y\rangle^{-q}(V_{g(v)}-V_{\mu_{0}})\|_{{\infty}}\lesssim
\delta_{0}b_{0}^{\frac{1+q}{2}}$$ which together with $\|\langle
y\rangle^{-q}(v-V_{\mu_{0}})\|_{{\infty}}\leq \delta_{0}b_{0}^{q/2}$
in the proposition implies (~\ref{eq:vh}).
\end{proof}
\section{A priori Estimates}\label{SEC:ApriEst}
% Proofread\\
In this section we assume that (~\ref{eq:MCF}) has a unique
solution, $u(x,t),$ $0\leq t\leq t_{*}$, such that
$v(y,\tau)=\lambda^{-\frac{2}{1-p}}(t)u(x,t)$, where $y=\lambda^{-1}
x$ and $\tau(t):=\int_{0}^{t}\lambda^{-2}(s)ds$, is in the
neighborhood $U_{\epsilon_{0}}$ determined in Proposition
~\ref{Prop:Splitting}. Then by Proposition ~\ref{Prop:Splitting}
there exist $C^{1}$ functions $a(\tau)$ and $b(\tau)$ such that
$v(y,\tau)$ can be represented as
\begin{equation}\label{eqn:split2}
v(y,\tau)=\lb\frac{1-p+by^{2}}{a+\frac{1}{2}}\rb^{\frac{1}{1-p}}+e^{\frac{ay^{2}}{4}}\xi(y,\tau)
\end{equation}
where $\xi(\cdot,\tau)\perp \phi_{0a},\phi_{2a}$ (see
(~\ref{eqn:split})). Now we set
$$-\lambda(t)\partial_{t}\lambda(t)=a(\tau(t)).$$

In the following we define some estimating functions to control
$\xi,$ $a$ and $b.$
\begin{equation}\label{eq:majorants}
\begin{array}{lll}
M_{1}(\tau):=\displaystyle\max_{s\leq \tau}\beta^{-3/2}(s)\|\langle
y\rangle^{-3}e^{\frac{ay^{2}}{4}}\xi(s)\|_{\infty},\\
M_{2}(\tau):=\displaystyle\max_{s\leq \tau}\beta^{-1}(s)\|\langle
y\rangle^{-2}e^{\frac{ay^{2}}{4}}\xi(s)\|_{\infty},\\
A(\tau):=\displaystyle\max_{s\leq \tau}\beta^{-2}(s)|a(s)-\frac{1}{2}+\frac{2}{1-p}b(s)|\\
B(\tau):=\displaystyle\max_{s\leq\tau}\beta^{-\frac{3}{2}}(s)|b(s)-\beta(s)|
\end{array}
\end{equation} moreover if $p<-1$
we define $$M_{q}(\tau):=\displaystyle\max_{s\leq
\tau}\beta^{-\frac{q}{2}}(\tau)\|\langle
y\rangle^{-q}e^{\frac{ay^{2}}{4}}\xi(s)\|_{\infty}$$ and recall that
$q:=\min\{\frac{4}{1-p},\frac{2(2-p)}{(p-1)^{2}},1\}.$ The function
$\beta(\tau)$ is defined as
\begin{equation}\label{FunBTau}
\beta(\tau):=\frac{1}{\frac{1}{b(0)}-\frac{4p}{(p-1)^{2}}\tau}.
\end{equation}

In the next we present a priori bounds on the fluctuation $\xi$
which are proved in later sections.
\begin{proposition}\label{Pro:Main2}
Suppose in (~\ref{eq:MCF}) the datum $u_{0}(x)$ satisfies all the
conditions in Theorem ~\ref{maintheorem}. Let the parameters
$a(\tau)$, $b(\tau)$ and the functions $v$ and $\xi$ the same as in
Proposition ~\ref{Prop:SplittingIC}. Then they satisfy the following
estimates: if
$$M_{1}(\tau)\leq \frac{1}{8}(1-p)^{\frac{1}{1-p}-\frac{3}{2}},\ A(\tau),\ B(\tau)\leq
\beta^{-1/4}(\tau)\ \text{and} \ v(y,\tau)\geq \frac{1}{2}
(\frac{1-p}{2})^{\frac{1}{1-p}}$$ in some interval $\tau\in [0,T]$
then for any $\tau\in [0,T]$ we have
\begin{enumerate}
\item[(1)]
if $0<p\leq -1$ then
\begin{equation}\label{MajorE}
B(\tau)\lesssim 1+M_1(\tau) A(\tau)+M_{1}^{2-p}(\tau)+A(\tau),
\end{equation}
\begin{equation}\label{EstA}
A(\tau)\lesssim
A(0)+1+\beta^{1/2}(0)[1+M_{1}^{2-p}(\tau)+A(\tau)M_{1}(\tau)],
\end{equation}
\begin{equation}\label{eq:M1}
M_{1}(\tau)\lesssim
M_{1}(0)+M_{1}(\tau)M_{2}(\tau)+\beta^{1/4}(0)[1+M_{1}^{2-p}(\tau)+A(\tau)M_{1}(\tau)],
\end{equation}
\begin{equation}\label{eq:M2}
M_{2}(\tau)\lesssim
M_{2}(0)+M_{1}(\tau)+M_{2}^{2}(\tau)+\beta^{1/4}(0)[1+M_{1}^{2-p}(\tau)+A(\tau)M_{1}(\tau)];
\end{equation}
\item[(2)]
if $p<-1$ then (~\ref{MajorE}) and (~\ref{EstA}) still hold, and
\begin{equation}\label{eq:hM1}
\begin{array}{lll}
M_{1}(\tau)&\lesssim&
M_{1}(0)+\beta^{1/4}(0)[1+M_{1}^{2-p}(\tau)+A(\tau)M_{1}(\tau)]+M_{1}(\tau)[M_{q}^{1-p}(\tau)+M_{q}(\tau)],
\end{array}
\end{equation}
\begin{equation}\label{eq:hM2}
\begin{array}{lll}
M_{2}(\tau)\lesssim
M_{2}(0)+M_{1}(\tau)+M_{2}(\tau)(M_{q}^{1-p}(\tau)+M_{q}(\tau))+\beta^{1/4}(0)[1+M_{1}^{2-p}(\tau)+M_{2}(\tau)+A(\tau)M_{1}(\tau)],
\end{array}
\end{equation}
\begin{equation}\label{eq:M3}
M_{q}(\tau)\lesssim
M_{q}(0)+M_{2}(\tau)+M_{q}^{2-p}(\tau)+M_{q}^{2}(\tau)+
\beta^{1/4}(0)[1+M_{q}(\tau)+M_{1}^{2-p}(\tau)+M_{1}(\tau)A(\tau)]
\end{equation}
\end{enumerate}
where the function $\beta$, the constant $q$ and the estimating functions
are defined above.
\end{proposition} We prove Equations (~\ref{MajorE}) and (~\ref{EstA})  in
Section ~\ref{SEC:EstB}, (~\ref{eq:M1}) and (~\ref{eq:hM1}) in Section ~\ref{SEC:EstM1}, (~\ref{eq:M2}) and
(~\ref{eq:hM2}) in Section ~\ref{SEC:EstM2},
(~\ref{eq:M3}) in Section ~\ref{SEC:EstM3}.
\section{The Lower Bound of $v$}\label{SEC:LowerBound}
% Proofread\\
In this section we prove a lower bound for $v$ defined in
(~\ref{eq:definev}). The main tool is a generalized form of maximum
principle.
\begin{lemma}\label{LM:max}
Suppose $u(y, \tau )$ is a smooth function having the following
properties: there exist smooth, bounded functions $a_{1}, a_{2}, d$
such that if $|y|\geq c( \tau )\geq 0$ then
$$u_{ \tau }-u_{yy}-[a_{1}(y, \tau )+d( \tau )y]u_{y}-a_{2}(y, \tau )u\leq 0;$$
$$\langle y\rangle^{-l}u(y, \tau )\in L^{\infty}  \ \text{with} \ l\geq 0;$$
$$u(y,0)\leq 0 \ \text{if}\ |y|\geq c( 0),\ \text{and}\ u(c( \tau ),
 \tau )\leq 0 \ \text{if}\  \tau \leq T;
$$ then we have $$u(y, \tau )\leq 0\ \text{if}\ |y|\geq c( \tau ),\  \tau \leq T.$$
\end{lemma}
\begin{proof}
In order to use the standard maximum principle we
first transform the function $u$. Define a new function $w$ by
\begin{equation}\label{Eq:defineV2}
\langle z\rangle^{l}w(z, \tau ):=u(y, \tau )
\end{equation} with $z:=y
e^{ \int_{0}^{\tau} d(s)ds}.$ Then $w$ is a smooth, bounded function
satisfying the inequality
$$w_{ \tau }-w_{zz}-a_{3}(z, \tau )w_{z}-a_{4}(z, \tau )w\leq 0$$
for some bounded, smooth functions $a_{3},\ a_{4};$ moreover
$$w(z,0)\leq 0 \ \text{if}\ |z|\geq c(0),\ \text{and}\ w(c( \tau ) e^{\int_{0}^{ \tau }d(s)ds},
 \tau )\leq 0 \ \text{if}\  \tau \leq T.$$ By the standard maximum principle
(see ~\cite{LSU}) we have $$w(z, \tau )\leq 0\ \text{if}\  \tau \leq
T,\ |z|\geq c( \tau )e^{\int_{0}^{ \tau }d(s)ds}.$$ By the definition of $w$ in (~\ref{Eq:defineV2}) we complete the proof.
\end{proof}
Recall the definition of the function $g(y,\beta)$ as
$$g(y,\beta):=\left\{
\begin{array}{lll}
(\frac{1-p}{2})^{\frac{1}{1-p}}\ \ \ \ \ \ \text{if}\ \beta y^{2}\leq 4(1-p)\\
(2(1-p))^{\frac{1}{1-p}}\ \text{if}\ \beta y^{2}>4(1-p)
\end{array}
\right.$$ The following lemma states an observation important for
our analysis.
\begin{lemma}\label{LM:UppLow} Let $v$ be as in (~\ref{eq:definev}).
Suppose that for time $\tau\leq \tau_{1}$, $M_{1}(\tau)\leq
\frac{1}{8}(1-p)^{\frac{1}{1-p}-\frac{3}{2}}$, $A(\tau), B(\tau)\leq
\beta^{-1/4}(\tau)$ and
\begin{equation}\label{eq:space}
\langle y\rangle^{-\frac{2}{1-p}}v(y,\tau)\in L^{\infty} \
\text{and}\ v(y,\tau) \geq c(\tau)
\end{equation} for some $c(\tau)> 0.$
Then we have
\begin{equation}\label{eq:comparison}
v(y,\tau)\geq (2c_{0}+\frac{2b_{0}}{1-p})^{\frac{1}{p-1}}g((2c_{0}+\frac{2b_{0}}{1-p})^{1/2}y,\beta(\tau)).
\end{equation}
\end{lemma}
\begin{proof}
By the scaling invariance of (~\ref{eq:MCF}), without losing
generality we assume that $2c_{0}+\frac{2b_{0}}{1-p}=1.$ By the
assumption on the datum we have that
\begin{equation}\label{eq:comparison0}
v(y,0)\geq g(y,\beta(0)).
\end{equation} Recall that $p<0$ and
$$v(y,\tau)=\lb\frac{1-p+b(\tau)y^{2}}{a(\tau)+\frac{1}{2}}\rb^{\frac{1}{1-p}}+e^{\frac{ay^{2}}{4}}\xi(y,\tau).$$
The assumption $M_{1}(\tau)\leq
\frac{1}{8}(1-p)^{\frac{1}{1-p}-\frac{3}{2}}$ yields
$|e^{\frac{a(\tau)y^{2}}{4}}\xi(y,\tau)|\leq
(\frac{1-p}{10})^{\frac{1}{1-p}}\beta^{3/2}(\tau)\langle
y\rangle^{3}$. Moreover by the assumption on $A$ we have $a(\tau)\in
[\frac{1}{4},\frac{3}{4}].$ Consequently
\begin{equation}
v(y,\tau)\geq
(\frac{1-p}{2})^{\frac{1}{1-p}} \ \text{when}\ \beta(\tau) y^{2}
\leq 4(1-p)
\end{equation}
and
\begin{equation}
v(y,\tau)\geq (2(1-p))^{\frac{1}{1-p}}\ \text{when}\ \beta(\tau) y^{2}
= 4(1-p)
\end{equation}
provided that $\beta(0)$ is sufficiently small.

For the region $\beta(\tau) y^{2} \geq 4(1-p)$, by the fact $p<0$ we
have that
\begin{equation}\label{eq:function}
H((2(1-p))^{\frac{1}{1-p}})\leq 0\ \text{and}\ H(v)=0
\end{equation} where the map $H(g)$
is defined as
$$H(g)=g_{\tau}-g_{yy}+g^{p}+ay\partial_{y}g-\frac{2a}{1-p}g.$$

Equations (~\ref{eq:comparison0})-(~\ref{eq:function}) and the
assumption (~\ref{eq:space}) enable us to use Lemma ~\ref{LM:max} on
the equation for $(2(1-p))^{\frac{1}{1-p}}-v$ to obtain
\begin{equation}\label{eq:lower2}
v(y,\tau)\geq g(y,\beta(\tau))\ \text{if}\ \beta(\tau)y^{2}\geq
4(1-p).
\end{equation} This together with the analysis on the region $\beta(\tau) y^{2}\leq
4(1-p)$ completes the proof.
\end{proof}
\section{Proof of Main Theorem ~\ref{maintheorem}}\label{SecMain}
% Proofread\\
In the next lemma we show that restriction (~\ref{eq:INI2}) on the
initial conditions involving two parameters can be rescaled into a
condition involving one parameter. Recall that
$q=\min\{\frac{4}{1-p},\ \frac{2(2-p)}{(1-p)^{2}},\ 1\}.$
\begin{lemma}\label{LM:rescale}
Let $u_{0}$ satisfy the condition (~\ref{eq:INI2}). Then there exist
some scalars $k_{0},\ \delta_{1},\ \beta>0$ such
that
\begin{equation}\label{eq:lowerbound}
\|\langle k_{0}x\rangle^{-n}[k_{0}^{\frac{2}{p-1}}u_{0}(k_{0}x)-
(\frac{1-p+\beta
x^{2}}{1-\frac{2}{1-p}\beta})^{\frac{1}{1-p}}]\|_{\infty}\leq
\delta_{1}\beta^{\frac{n}{2}},
\end{equation}
$$k_{0}^{\frac{2}{p-1}}u_{0}(k_{0}x)\geq g(x,\beta)$$ for $n=2,3$, and
if $p<-1$ (~\ref{eq:lowerbound}) with $n=q.$
\end{lemma}
\begin{proof}
Define $k_{0}:=(2 c_{0}+\frac{2}{1-p}b_{0})^{-1/2}$,
$\delta_{1}:=\delta_{0}k_{0}^{\frac{2}{p-1}}$ and
$\beta:=b_{0}k_{0}^{2}$. It is straightforward to verify that the
function $k_{0}^{\frac{2}{p-1}}u_{0}(k_{0}x)$ has all the properties
above.
\end{proof}
By this lemma in what follows we only study the case
\begin{equation}\label{eq:INI}
\begin{array}{rrr}
\|\langle
x\rangle^{-n}[u_{0}(x)-(\frac{1-p+b_{0}x^{2}}{1-\frac{2}{1-p}b_{0}})^{\frac{1}{1-p}}]\|_{\infty}\leq
\delta_{0}b_{0}^{\frac{n}{2}},
\end{array}
\end{equation} $$u_{0}(x)\geq g(x,b_{0})$$ for $n=2,3,$ and
if $p<-1$, (~\ref{eq:INI}) with $n=q$.

By Proposition ~\ref{THM:WellPose}, there exists $\infty\ge t_{*}>0$
such that Equation (~\ref{eq:MCF}) has a unique solution $u(x,t)$
for $0\leq t< t_*$, but has no solutions on a larger time interval.
Moreover, if $t_* < \infty$, then $\|
\frac{1}{u(\cdot,t)}\|_{\infty} \rightarrow \infty$ as $t\rightarrow
t_*$. Recall that the solutions $u(x,t)$, $v(y, \tau)$ and
$w(y,\tau)$ and the corresponding initial conditions are related by
the scaling and gauge transformations (see (~\ref{eq:definev}) and
(~\ref{eq:definew})). Take $\lambda(0)=1$. Then we have that
$u_{0}(x)=v_{0}(y)$.

Choose $b_0$ so that $\delta_{0}b_{0}^{3/2}\leq
\frac{1}{2}\epsilon_{0}$ with $\delta_{0}$ the same as in
\eqref{eq:INI} and with $\epsilon_{0}$ given in Proposition
~\ref{Prop:Splitting}. Then $v_{0}\in U_{\frac{1}{2}\epsilon_{0}}$,
by the condition \eqref{eq:INI} with $n=3$ on the initial
conditions.
%Take $b_{0}=\epsilon_{0}$ as given in Proposition ~\ref{Prop:Splitting}.
By continuity there is a (maximal) time $t_\#\leq t_{*}$ such that
$v \in U_{\epsilon_0}$ for $t< t_\#$. For this time interval
Propositions \ref{Prop:Splitting} and \ref{Prop:SplittingIC} hold
for $v$ and in particular we have the splitting \eqref{eqn:split2}.
Recall that we assume $a=2c-\frac{1}{2}$ in the decomposition
\eqref{eqn:split2}. In particular, this implies that the initial
condition can be written in the form
\begin{equation} \label{eqn:splitic}
v_0(y) =V_{a(0) b(0)}(y) + e^{\frac{a(0)y^{2}}{4}}\xi_0(y),
\end{equation}
where $(a(0), b(0))=g(v_0)$ and $\xi_0\perp e^{-\frac{a(0)}{4} y^2},
(1-a(0) y^2)e^{-\frac{a(0)}{4} y^2}$. Moreover $M_{n}(0)\lesssim
\delta_{0},\ n=1,2$ and especially $q$ for $p<-1.$

Before starting proof, we note that the definition of $g$ in Lemma
~\ref{LM:UppLow} implies
$$g(y,\tau)\geq (\frac{1-p}{2})^{\frac{1}{1-p}}$$ for any
time $\tau.$

By the conditions of Theorem ~\ref{maintheorem} on the datum, we
have that $b(0)$ is small, $M_{1}(0),\ M_{2}(0)\lesssim \delta_{0}$ for
some small $\delta_{0}$, $A(0)$ and $B(0)$ are bounded, and
$$u_{0}\in \langle x\rangle^{\frac{2}{1-p}}L^{\infty}\ \text{and}\
u_{0}\geq (\frac{1-p}{2})^{\frac{1}{1-p}}.$$ Then by the local well posedness Theorem ~\ref{THM:WellPose} and the splitting Proposition ~\ref{Prop:SplittingIC}
there exists a time interval $[0,T]$ such that for any $\tau\in
[0,T]$,
$$M_{1}(\tau)\leq \frac{1}{8}(1-p)^{\frac{1}{1-p}-\frac{3}{2}},\ B(\tau),\ A(\tau)\leq
\beta^{-1/4}(\tau)$$ and $v(y,\tau)\geq
\frac{1}{2}(\frac{1-p}{2})^{\frac{1}{1-p}}.$ Thus we have Equations
(~\ref{MajorE})-(~\ref{eq:M3}) which together with the fact that $\beta(0)$ is small imply that
$$M_{1}(\tau), \ M_{2}(\tau), M_{q}(\tau)\lesssim \delta_{0},\ B(\tau),\ A(\tau)\lesssim 1 \ll
\beta^{-1/4}(\tau)$$ if $\tau\in [0,T].$ In the next we only show
how to use (~\ref{EstA})-(~\ref{eq:M2}) to get the estimate for $A$,
$B$, $M_{1}$ and $M_{2}$, the proof of $M_{q}$ is similar. Indeed,
since $M_{1}(\tau)\leq 1$, we can solve \eqref{EstA} for $A(\tau)$.
We substitute the result into Equations \eqref{eq:M1} -
\eqref{eq:M2}, and substitute the estimate for $M_{1}(\tau)$ into
the right hand side of (~\ref{eq:M2}), to obtain inequalities
involving only the estimating functions $M_{1}(\tau)$ and
$M_{2}(\tau)$. Consider the resulting inequality for $M_{2}(\tau)$.
The only terms on the r.h.s., which do not contain $\beta(0)$ to a
power at least $1/4$ as a factor, are $M_{2}^{2}(\tau)$ and
$M_{1}(\tau)M_{2}(\tau)$. Hence for $M_{2}(0)\ll 1$ this inequality
implies that $M_{2}(\tau) \lesssim M_{2}(0)+M_{1}(0)  +
\beta^{\frac{1}{4}}(0)$. Substituting this result into the
inequality for $M_{1}(\tau)$ we obtain that $M_{1}(\tau) \lesssim
M_{1}(0)  + \beta^{\frac{1}{4}}(0)$ as well. The last two
inequalities together with \eqref{MajorE} and \eqref{EstA} imply the
desired estimates on $A(\tau)$ and $B(\tau)$.

Moreover by Lemma ~\ref{LM:UppLow} we obtain
$$v(\cdot,\tau)\in \langle y\rangle^{\frac{2}{1-p}} L^{\infty}\ \text{and}\ v(y,
\tau)\geq g(y,\tau)\geq (\frac{1-p}{2})^{\frac{1}{1-p}}.$$ By using
the procedure above recursively we have that
$$M_{1}(\tau), \ M_{2}(\tau), M_{q}(\tau)\lesssim \delta_{0},\ B(\tau),\ A(\tau)\lesssim 1$$ for any
$\tau.$

By the definitions of $A(\tau)$ and $B(\tau)$ in
(~\ref{eq:majorants}) and the facts that $A(\tau), B(\tau)\lesssim
1$ proved above, we have
\begin{equation}\label{EstABtau}
a-\frac{1}{2}=-\frac{2b}{1-p}+O(\beta^{2}),\ b(\tau)=\beta(\tau)+O(
\beta^{3/2}).
\end{equation}
Hence $a-\frac{1}{2}=O(\beta(\tau)),$ where, recall the definition
of $\beta(\tau)$ from (~\ref{FunBTau}).

The facts $\lambda(0)=1$ and $a=-\lambda\frac{d}{dt}\lambda$ after
Equation (~\ref{eqn:BVNLH}) yield
\begin{equation}
\lambda(t)=[1-2\int_{0}^{t}a(\tau(s))ds]^{1/2}.
\end{equation} Since $|a(\tau(t))-\frac{1}{2}|=O(b(\tau(t)))$ there exists a time $t^{*}$ such
that $1=2\int_{0}^{t^*}a(\tau(s))ds$, i.e. $\lambda(t)\rightarrow 0$
as $t\rightarrow t^{*}$ with the latter defined in (~\ref{eq:INI}).
For any time $t\leq t^*$ (~\ref{eq:MCF}) is well posed and can be
split by the fact that $M_{1}(\tau(t)),\ M_{2}(\tau(t)),\ M_{q}(\tau(t))\lesssim
\delta_{0} $ and $ A(\tau(t))$ and $B(\tau(t))$ are bounded.
Obviously $t^{*}\leq t_{*}$. Moreover by the definition of $\tau$
and the property of $a$ we have that $\tau(t)\rightarrow \infty$ as
$t\rightarrow t^{*}.$ Equation (~\ref{EstABtau}) implies
$b(\tau(t))\rightarrow 0$ and $a(\tau(t))\rightarrow \frac{1}{2}$ as
$t\rightarrow t^*,$ where, recall the definition of $\beta(\tau)$
from (~\ref{FunBTau}). By the analysis above and the definitions of
$a$, $\tau$ and $\beta$ we have
$$\lambda(t)=(t^*-t)^{1/2}(1+o(1)),\ \tau(t)=-ln|t^*-t|(1+o(1))$$ with $o(1)\rightarrow 0$ as $t\rightarrow t^*$
and consequently
$$
\beta(\tau(t))=\frac{(p-1)^{2}}{4p|ln|t^*-t||}(1+O(\frac{1}{|ln|t^*-t||^{1/2}}))$$
By (~\ref{EstABtau}) we have $$
b(\tau(t))=\frac{(p-1)^{2}}{4p|ln|t^*-t||}(1+O(\frac{1}{|ln|t^*-t||^{1/2}}))$$
and $$
a(\tau(t))=\frac{1}{2}+\frac{1-p}{2p|ln|t^*-t||}(1+O(\frac{1}{ln|t^*-t|})).
$$ If we define $2c(t)=a(\tau(t))+\frac{1}{2}$ in Theorem ~\ref{maintheorem}, then $c(t)$ has the same
estimate as (~\ref{eq:para}).

By the fact $M_{1}\ll 1$ we have $$\|\langle y\rangle^{-3}e^{\frac{ay^{2}}{4}}\xi(\cdot,\tau(t))\|_{\infty}\ll \beta^{3/2}(\tau(t))$$thus
$$u(0,t)=\lambda^{\frac{2}{1-p}}(t)[V_{a,b}+e^{\frac{ay^{2}}{4}}\xi(y,\tau)]|_{y=0}
\leq
\lambda^{\frac{2}{1-p}}(t)[(\frac{1-p}{a(\tau(t))+\frac{1}{2}})^{\frac{1}{1-p}}+C\beta^{3/2}(\tau(t))]$$
for some constant $C>0,$ hence $\frac{1}{u(0,t)}\rightarrow \infty$
as $t\rightarrow t^*.$ Moreover by the facts $v(y,\tau)\geq
(\frac{1-p}{2})^{\frac{1}{1-p}}$ proved above and
$u(x,t)=\lambda^{\frac{2}{1-p}}(t)v(y,\tau)$ we have that for any
$t<t^*$ $$\| \frac{1}{u(\cdot,t)}\|_{\infty}\leq
\lambda^{\frac{2}{p-1}}(t)(\frac{2}{1-p})^{\frac{1}{1-p}}<\infty.$$
Those facts together with the estimates of $\lambda(t),\ b(\tau(t))$
and $a(\tau(t))$ imply that $\frac{1}{u(\cdot,t)}$ collapses at time
$t^{*}$ and $t^*=t_*$. Thus we finish proving the first argument of
Theorem ~\ref{maintheorem}.

Collecting the facts above we complete the proof.
\section{Parametrization and the Linearized
Operator}\label{SEC:ParaMetrization}
Substitute (~\ref{eqn:split2}) into (~\ref{eqn:BVNLH})
we have the following equation for $\xi$
\begin{equation}\label{eq:xi}
\frac{d}{d\tau}\xi(y,\tau)=-L(a,b)\xi+F(a,b)+N(a,b,\xi)
\end{equation}
where the linear operator $L(a,b)$ is defined as
$$L(a,b):=-\partial_{y}^{2}+\frac{a^{2}+a_{\tau}}{4}y^{2}-\frac{a}{2}-\frac{2a}{1-p}+p\frac{\frac{1}{2}+a}{1-p+b y^{2}},$$ and the function $F(a,b)$ is defined as
\begin{equation}\label{eq:scource}
F(a,b):=(\frac{1-p+by^{2}}{a+\frac{1}{2}})^{\frac{1}{1-p}}\frac{1}{1-p}e^{-\frac{ay^{2}}{4}}[\Gamma_{1}+\Gamma_{2}\frac{y^{2}}{1-p+by^2}+F_{1}]
\end{equation}
with $$
\Gamma_{1}:=\frac{a_{\tau}}{a+\frac{1}{2}}+a-\frac{1}{2}+\frac{2b}{1-p},$$
$$
\Gamma_{2}:=-b_{\tau}-b(a-\frac{1}{2}+\frac{2b}{1-p})+\frac{4p}{(p-1)^{2}}b^{2},
$$
\begin{equation}\label{eq:MainRemain}
F_{1}:=\frac{p}{(1-p)^{2}}\frac{4b^{3}y^{4}}{(1-p+by^{2})^{2}};
\end{equation}
and the nonlinear term $N$ is
$$
N(a,b,\xi):=-v^{p}e^{-\frac{ay^{2}}{4}}+(\frac{a+\frac{1}{2}}{1-p+by^{2}})^{-\frac{p}{1-p}}e^{-\frac{ay^{2}}{4}}+p\frac{a+\frac{1}{2}}{1-p+by^{2}}\xi,
$$
where, recall the definition of $v$ in (~\ref{eq:definev}).

In the next we derive various estimates for $\Gamma_{1},$ $\Gamma_{2}$ and $N(a,b,\xi)$.
\begin{lemma}\label{LM:ESTnonlin}
Suppose that $A(\tau),B(\tau)\leq \beta^{-1/4}(\tau)$,
$v(y,\tau)\geq \frac{1}{2}(\frac{1-p}{2})^{\frac{1}{1-p}}$ and $b(0)=\beta(0)$ is small. Then for
$0> p\geq -1$ we have
\begin{equation}\label{eq:ESTnonlin}
|N(a,b,\xi)|\lesssim e^{-\frac{ay^{2}}{4}}(\frac{1}{1+\beta
y^{2}})^{\frac{2}{1-p}}|e^{\frac{ay^{2}}{4}}\xi|^{2};
\end{equation} and for $p<-1$
\begin{equation}\label{eq:pESTnonlin}
|N(a,b,\xi)|\lesssim e^{-\frac{ay^{2}}{4}}[(\frac{1}{1+\beta
y^{2}})^{\frac{2-p}{1-p}}|e^{\frac{ay^{2}}{4}}\xi|^{2-p}+(\frac{1}{1+\beta
y^{2}})^{\frac{2}{1-p}}|e^{\frac{ay^{2}}{4}}\xi|^{2}];
\end{equation}
\end{lemma}
\begin{proof}
By direct computation and the assumption on $v$ we have
$$
\begin{array}{lll}
|N(a,b,\xi)|&\lesssim & e^{-\frac{ay^{2}}{4}}v^{p}
[|(1+\phi)^{-p}-1+p\phi|+p|\phi||(1+\phi)^{-p}-1|]\\
& \lesssim & e^{-\frac{ay^{2}}{4}}
[|(1+\phi)^{-p}-1+p\phi|+p|\phi||(1+\phi)^{-p}-1|]
\end{array}
$$ where the function $\phi$ is defined as
$\phi:={e^{\frac{ay^{2}}{4}}\xi}{V^{-1}_{a,b}}$ with
$V_{a,b}:=(\frac{1-p+b y^{2}}{\frac{1}{2}+a})^{\frac{1}{1-p}}.$
\begin{enumerate}
\item[(A)]
When $|\phi|\geq 1/2$ we have that if $0>p\geq -1$ then $$1,
|(1+\phi)^{-p}| \lesssim \phi^{2},\ \text{thus}\
|N(a,b,\xi)|\lesssim e^{-\frac{ay^{2}}{4}}\phi^{2};$$ if $p<-1$ then
$$1,\ |(1+\phi)^{-p}|,\ |\phi|\lesssim |\phi|^{2-p},\ \text{thus}\ |N(a,b,\xi)|\lesssim e^{-\frac{ay^{2}}{4}}|\phi|^{2-p};$$
\item[(B)]
when $|\phi|< 1/2$ by the remainder estimate of the Tylor
expansion we have
$$|(1+\phi)^{-p}-1+p\phi|, \ |\phi||(1+\phi)^{-p}-1|\lesssim
\phi^{2}. $$
\end{enumerate}
Collecting the estimate above we have that if $0>p\geq -1$ then
$$|N(a,b,\xi)|\lesssim e^{-\frac{ay^{2}}{4}} |\phi|^{2};$$ if $p<-1$
then $$|N(a,b,\xi)|\lesssim e^{-\frac{ay^{2}}{4}}
(|\phi|^{2}+|\phi|^{2-p}).$$

Moreover the assumptions $A(\tau),B(\tau)\leq \beta^{-1/4}(\tau)$
imply that $b=\beta+o(\beta)$ and $a=\frac{1}{2}+O(\beta)$, thus
$V_{a,b}\gtrsim ({1+\beta(\tau)y^{2}})^{\frac{1}{1-p}}$, i.e.
$$|\phi|\lesssim (\frac{1}{1+\beta y^{2}})^{\frac{1}{1-p}}e^{\frac{a
y^{2}}{4}}\xi$$ which together with the definition of $\phi$ and the
estimates of $N$ above implies (~\ref{eq:ESTnonlin}).
\end{proof}
In the following proposition we establish
some estimates for $\Gamma_{1}$ and $\Gamma_{2}$.
\begin{proposition}
Suppose that $M_{1}(\tau)\leq \frac{1}{8}(1-p)^{\frac{1}{1-p}-\frac{3}{2}},\
A(\tau),\ B(\tau)\leq \beta^{-1/4}(\tau)$ and $v(y,\tau)\geq
\frac{1}{2}(\frac{1-p}{2})^{\frac{1}{1-p}}$ and $b(0)=\beta(0)$ is small. Then we have
\begin{equation}\label{eq:estI1I2}
|\Gamma_{1}|,\ |\Gamma_{2}| \lesssim
\beta^{5/2}(1+M^{2-p}_{1}+AM_{1}).
\end{equation}
\end{proposition}
\begin{proof}
Taking inner products on (~\ref{eq:xi}) with the functions
$\phi_{0,a}:=(\frac{\alpha}{2\pi})^{\frac{1}{4}}e^{-\frac{ay^{2}}{4}}$ and $\phi_{2,a}:=(\frac{\alpha}{8\pi})^{\frac{1}{4}}(ay^{2}-1)e^{-\frac{ay^{2}}{4}}$ and
using the facts $\xi\perp \phi_{0,a},\ \phi_{2,a}$ in (~\ref{eqn:split2}), we have
\begin{equation}\label{eq:Projection}
|\langle F(a,b),\phi_{0,a}\rangle|\leq  G_{1},\ \
|\langle F(a,b),\phi_{2,a}\rangle|\leq  G_{2}
\end{equation} where the functions $G_{1},\ G_{2}$ are defined as
$$G_{1}:=|\langle by^{2}\xi,\phi_{0,a}\rangle|+|\langle a_{\tau}y^{2}\xi,\phi_{0,a}\rangle|+|\langle
N(a,b,\xi),\phi_{0,a}\rangle|,$$ and $$G_{2}:=|\langle
by^{2}\xi,\phi_{2,a}\rangle|+|\langle
a_{\tau}y^{2}\xi,\phi_{2,a}\rangle|+|\langle
N(a,b,\xi),\phi_{2,a}\rangle|.$$ The estimate
of $N$ in (~\ref{eq:ESTnonlin}), the assumption on $B(\tau)$ and the
definition of $M_{1}$ yield
$$G_{1},G_{2}\lesssim
\beta^{5/2}[M_{1} +M^{2}_{1}+M_{1}^{2-p}]+|a_{\tau}|\beta^{3/2}M_{1}.$$

Now we study the right hand sides of (~\ref{eq:Projection}). We
rewrite the function $F(a,b)$ in (~\ref{eq:scource}) as
$$F(a,b)=\chi(a,b)[\Gamma_{1}+\Gamma_{2}\frac{1}{a[1-p+by^{2}]}+\Gamma_{2}\frac{ay^{2}-1}{a[1-p+by^{2}]}+F_{1}].$$
where, recall the definition of $F_{1}$ in (~\ref{eq:MainRemain}),
and the term $\chi(a,b)$ is defined as
$$\chi(a,b):=(\frac{1-p+by^{2}}{a+\frac{1}{2}})^{\frac{1}{1-p}}\frac{1}{1-p}e^{-\frac{ay^{2}}{4}}.$$
In the $L^{\infty}$ space we expand $F(a,b)$ as
$$
F(a,b)=(\frac{1}{1-p})^{\frac{p}{p-1}}(a+\frac{1}{2})^{\frac{1}{p-1}}e^{-\frac{a
y^{2}}{4}}[\Gamma_{1}+\frac{1}{a(1-p)}\Gamma_{2}+\frac{ay^{2}-1}{a(1-p)}\Gamma_{2}]+O(b^{3}+|b|[|\Gamma_{1}|+|\Gamma_{2}|])
$$ where we use the observation $\|\chi(a,b)F_{1}\|_{\infty}=O(b^{3})$ and the fact $a\geq \frac{1}{4}$ implied by the assumption on $A.$

By using the fact that $\phi_{0,a}\perp
\phi_{2,a}$ we have
$$|\langle F(a,b),\phi_{0,a}\rangle|\gtrsim
|\Gamma_{1}+\frac{1}{a(1-p)}\Gamma_{2}|-|b|(|\Gamma_{1}|+|\Gamma_{2}|)-|b|^{3}$$
and
$$|\langle F(a,b),\phi_{2,a}\rangle|\gtrsim
|\Gamma_{2}|-|b|(|\Gamma_{1}|+|\Gamma_{2}|)-b^{3},$$ which together
with (~\ref{eq:Projection}) and the assumption on $B$ implies that
\begin{equation}\label{eq:I1i2}
|\Gamma_{1}|,\ |\Gamma_{2}|\lesssim
\beta^{5/2}(1+M_{1}^{2}+M_{1}^{2-p})+|a_{\tau}|\beta^{2}M_{1}.
\end{equation}
By the definitions of $\Gamma_{1}$ and $A$, $a_{\tau}$ has the bound
$$|a_{\tau}|\lesssim |\Gamma_{1}|+\beta^{3/2}(\tau)A(\tau).$$ Consequently after using $M_{1}(\tau)\leq
\frac{1}{8}(1-p)^{\frac{1}{1-p}-\frac{3}{2}}$ and some manipulation on (~\ref{eq:I1i2}) we obtain
$$|\Gamma_{1}|,\ |\Gamma_{2}|\lesssim
\beta^{5/2}(1+M_{1}^{2-p}+\beta^{1/2}AM_{1}).$$

The proof is complete.
\end{proof}
Recall that $q:=\min\{\frac{4}{1-p},\frac{2(2-p)}{(p-1)^2},1\}.$ The
estimates on $\Gamma_{1}$ and $\Gamma_{2}$ imply
\begin{corollary}
\begin{equation}\label{eq:estSource}
\|\langle y\rangle^{-n}e^{\frac{ay^{2}}{4}}F\|_{\infty} \lesssim
\beta^{\frac{n}{2}+1}[1+M^{2-p}_{1}+AM_{1}]
\end{equation} with $n=2,3$ and especially $n=q$ if $p<-1.$
\end{corollary}
\begin{proof}
Recall that $p<0$ and the definition of $F$ in Equation
(~\ref{eq:scource}) as
$$e^{\frac{ay^{2}}{4}}F(a,b)=(\frac{1-p+by^{2}}{a+\frac{1}{2}})^{\frac{1}{1-p}}\frac{1}{1-p}[\Gamma_{1}+\Gamma_{2}\frac{y^{2}}{1-p+by^2}+F_{1}].
$$ The estimates on $\Gamma_{1}$ and $\Gamma_{2}$ in (~\ref{eq:estI1I2})
implies that
$$
\begin{array}{lll}
\|\langle
y\rangle^{-n}(\frac{1-p+by^{2}}{a+\frac{1}{2}})^{\frac{1}{1-p}}[\Gamma_{1}+\Gamma_{2}\frac{y^{2}}{1-p+by^2}]\|_{\infty}
\lesssim  |\Gamma_{1}|+|\Gamma_{2}| \lesssim
\beta^{5/2}(1+M^{2-p}_{1}+AM_{1})
\end{array}
$$ with $n=2,3.$ Moreover for $p<-1$ we have that $1\geq q>\frac{2}{1-p}$ hence
$$\frac{|y|^{2-q}}{(1-p+by^{2})^{1+\frac{1}{p-1}}}\leq \frac{|y|^{2-q}}{(1-p+by^{2})^{\frac{2-q}{2}}}\leq
b^{\frac{q-2}{2}}$$ consequently
$$\|\langle y\rangle^{-q}(\frac{1-p+by^{2}}{a+\frac{1}{2}})^{\frac{1}{1-p}}[\Gamma_{1}+\Gamma_{2}\frac{y^{2}}{1-p+by^2}]\|_{\infty}\lesssim \beta^{\frac{q}{2}+1}(1+M^{2-p}_{1}+AM_{1})$$

For the term $F_{1}$ by straightforward computation we have
$$\|\langle y\rangle^{-2}F_{1}(\frac{1-p+by^{2}}{a+\frac{1}{2}})^{\frac{1}{1-p}}\|_{\infty}\lesssim
\beta^{2},\ \ \|\langle
y\rangle^{-3}F_{1}(\frac{1-p+by^{2}}{a+\frac{1}{2}})^{\frac{1}{1-p}}\|_{\infty}\lesssim
\beta^{5/2}(\tau),$$ using the fact $1\geq q>\frac{2}{1-p}$ for
$p<-1$ again
$$\|\langle y\rangle^{-q}F_{1}(\frac{1-p+by^{2}}{a+\frac{1}{2}})^{\frac{1}{1-p}}\|_{\infty}\lesssim
\beta^{1+\frac{q}{2}}.$$

Collecting the estimates above we complete the proof.
\end{proof}
\section{Proof of Estimates
\eqref{MajorE}-\eqref{EstA}}\label{SEC:EstB}
\begin{lemma}
If $B(\tau)\le \beta^{-1/4}(\tau)$ for $\tau\in[0,T]$, then Equation (~\ref{MajorE}) holds.
\end{lemma}
\begin{proof}
We rewrite the estimate of $\Gamma_{2}$ in Equation
\eqref{eq:estI1I2} as
\begin{equation*}
|b_\tau-\frac{4p}{(p-1)^{2}} b^2|\lesssim b| \frac{1}{2}-a-\frac{2
b}{p-1} |+\beta^{5/2}(1+M^{2-p}_{1}+AM_{1}).
\end{equation*}
The first term on the right hand side is bounded by $b \beta^{2}
A\lesssim \beta^{5/2}A$ by the definition of $A$, consequently
$$|b_{\tau}-\frac{4p}{(p-1)^{2}}b^{2}|\lesssim \beta^{5/2}(1+M_{1}A+M_{2}^{2-p}+M_{2}^{2}).$$
We divide both sides by $b^2$ and use the
inequality $b\lesssim \beta$ implied by the assumption on $B$ to
obtain the estimate
\begin{equation}
\left| \partial_\tau\frac{1}{b}+\frac{4p}{(p-1)^{2}} \right|\lesssim
\beta^{1/2} (1+ M_1 A+M_{1}^{2-p}+A). \label{est:InversebDE}
\end{equation}
Since $\beta$ is a solution to $\partial_\tau
\beta^{-1}+4p(p-1)^{-2}=0$, $\beta(0)=b(0)$, Equation
\eqref{est:InversebDE} implies that
\begin{equation*}
|\partial_\tau( \frac{1}{b}-\frac{1}{\beta})|\lesssim
\beta^{1/2}(1+M_1 A+M_{1}^{2-p}+A).
\end{equation*}
Integrating this equation over $[0,\tau]$, multiplying the result by
$\beta^{\frac{1}{2}}$ and using that $b\lesssim \beta$ give the
estimate
\begin{equation*}
\beta^{-\frac{3}{2}}|\beta-b|\lesssim \beta^{\frac{1}{2}}\int_0^\tau
\beta^{1/2}(s) (1+M_1 A+M_{1}^{2-p}+A) ds.
\end{equation*}
By the definitions of $\beta$ we have
$\beta^{1/2}(\tau)\int_{0}^{\tau}\beta^{1/2}(s)ds \lesssim 1$ which
together with the definition of $B$ gives (~\ref{MajorE}).
\end{proof}
\begin{lemma}
If $A(\tau),\ B(\tau)\le \beta^{-1/4}(\tau)$ for $\tau\in[0,T]$, then Equation (~\ref{EstA}) holds.
\end{lemma}
\begin{proof}
Define the quantity
$\Gamma:=a-\frac{1}{2}+\frac{2}{1-p} b$.  We prove the proposition
by integrating a differential inequality. Differentiating $\Gamma$
with respect to $\tau$ and substituting for $b_\tau$ and $a_\tau$ in
Equation \eqref{eq:estI1I2} we obtain
\begin{equation*}
\partial_\tau \Gamma+ ( a+\frac{1}{2}-\frac{2}{p-1} b
)\Gamma=-\frac{8 p}{(p-1)^3} b^2+{\cal R}_b.
\end{equation*} where ${\cal R}_{b}$ has the bound $$|\mathcal{R}_{b}|\leq
\beta^{5/2}(1+M_{1}^{2-p}+AM_{1}).$$ Let $\mu=\exp{\int_{0}^{\tau}
a(s)+\frac{1}{2}-\frac{2}{p-1} b(s) ds}$. Then the above equation
implies that
\begin{equation*}
\mu\Gamma-\Gamma(0)=-\frac{8p}{(p-1)^3}\int_0^\tau \mu b^2\, ds+\int_0^\tau \mu
{\cal R}_b\, ds.
\end{equation*}
We now use the inequality $b\lesssim \beta$ and the estimate of
${\cal R}_b$ to simplify the bound of $\Gamma$ as
\begin{equation*}
|\Gamma|\lesssim \mu^{-1}\Gamma(0)+\mu^{-1}\int_0^\tau \mu \beta^2\,
ds+\mu^{-1}\int_0^\tau \mu \beta^{5/2}(1+M^{2-p}_{1}+AM_{1}) ds.
\end{equation*}
For our purpose, it is sufficient to use the less sharp inequality
\begin{equation*}
|\Gamma|\lesssim\mu^{-1}\Gamma(0)+\mu^{-1}\int_0^\tau \mu
\beta^{2}\, ds ( 1+\beta^{1/2}(0)[M^{2-p}_{1}+AM_{1}] ).
\end{equation*} The assumption that $A(\tau),\ B(\tau)\leq
\beta^{-1/4}(\tau)$ implies that $a+\frac{1}{2}-\frac{2}{p-1}b\geq
\frac{1}{2}.$ Thus it is not difficult to show that
$\beta^{-2}\mu^{-1}\Gamma(0)\leq A(0);$ and by the slow decay of
$\beta$ we have that if $s\leq \tau$ then
$e^{-\frac{\tau-s}{4}}\beta^{2}(s)\lesssim \beta^{2}(\tau)$,
consequently $$
\beta^{-2}(\tau)\mu^{-1}(\tau)\int_0^\tau \mu(s) \beta^{2}(s)\,
ds\leq
\beta^{-2}(\tau)\int_{0}^{\tau}e^{-\frac{\tau-s}{2}}\beta^{2}(s)ds
\lesssim \int_{0}^{\tau}e^{-\frac{\tau-s}{4}}ds \lesssim 1.
$$ Collecting the estimates above we have $$\beta^{-2}|\Gamma|\lesssim 1+A(0)+\beta^{1/2}(0)(1+M_{1}^{2-p}+AM_{1})$$ which together with the definition of $A$ implies (~\ref{EstA}). Thus the proof is complete.
\end{proof}
\section{Rescaling of Fluctuations on a Fixed Time
Interval}\label{SEC:Rescale}
We return to our key equation (~\ref{eq:xi}). In this section we
re-parameterize the unknown function $\xi(y,\tau)$ in such a way
that the $y^{2}$-term in the linear part of the new equation has a
time-independent coefficient (cf ~\cite{DGSW}).

Let $t(\tau)$ be the inverse function to $\tau(t)$, where
$\tau(t)=\int_{0}^{t}\lambda^{-2}(s)ds$ for any $\tau\geq 0.$ Pick
$T>0$ and approximate $\lambda(t(\tau))$ on the interval $0\leq
\tau\leq T$ by the new trajectory, $\lambda_{1}(t(\tau)),$ tangent
to $\lambda(t(\tau))$ at the point $\tau=T:$
$\lambda_{1}(t(T))=\lambda(t(T)),$ and
$\alpha:=-\lambda_{1}(t(\tau))\partial_{t}\lambda_{1}(t(\tau))=a(T)$
where, recall
$a(\tau):=-\lambda(t(\tau))\partial_{t}\lambda(t(\tau)).$ Now we
introduce the new independent variables $z$ and $\sigma$ as
$z(x,t):=\lambda^{-1}_{1}(t)x$ and
$\sigma(t):=\int_{0}^{t}\lambda^{-2}_{1}(s)ds$ and the new unknown
function $\eta(z,\sigma)$ as
\begin{equation}\label{NewFun}
\lambda_{1}^{\frac{2}{1-p}}(t)e^{\frac{\alpha
z^{2}}{4}}\eta(z,\sigma):=\lambda^{\frac{2}{1-p}}(t)e^{\frac{ay^{2}}{4}}\xi(y,\tau).
\end{equation}
In this relation one has to think of the variables $z$ and $y$,
$\sigma$, $\tau$ and $t$ as related by
$z=\frac{\lambda(t)}{\lambda_{1}(t)}y,$
$\sigma(t):=\int_{0}^{t}\lambda_{1}^{-2}(s)ds$ and
$\tau=\int_{0}^{t}\lambda^{-2}(s)ds,$ and moreover
$a(\tau)=-\lambda(t(\tau))\partial_{t}\lambda(t(\tau))$ and
$\alpha=a(T).$

For any $\tau=\int_{0}^{t(\tau)}\lambda^{-2}(s)ds$ with $t(\tau)\leq
t(T)$ (or equivalently $\tau\leq T$) we define a new function
$\sigma(\tau):=\int_{0}^{t(\tau)}\lambda_{1}^{2}(s)ds.$ Observe the
function $\sigma$ is invertible, we denote by $\tau(\sigma)$ as its
inverse. Especially we define
\begin{equation}\label{T2}
S:=\int_{0}^{t(T)}\lambda_{1}^{-2}(s)ds.
\end{equation}
The new function $\eta$ satisfies the equation
\begin{equation}\label{eq:eta}
\frac{d}{d\sigma}\eta(\sigma)=-\mathcal{L}_{\alpha,\beta}\eta(\sigma)+\mathcal{W}\eta(\sigma)+\mathcal{F}(a,b)(\sigma)+\mathcal{N}_{1}(a,b,\alpha,\eta)(\sigma),
\end{equation}
with the operators
$$\mathcal{L}_{\alpha,\beta}:=L_{\alpha}+V,$$
$$L_{\alpha}:=-\frac{d^{2}}{dz^{2}}+\frac{\alpha^{2}}{4}z^{2}-\frac{1}{2}\alpha-\frac{2\alpha}{1-p},$$
$$V:=\frac{2 p \alpha}{1-p+\beta(\sigma) z^{2}},$$
$$\mathcal{W}:=p\frac{\lambda^{2}}{\lambda_{1}^{2}}\frac{a+\frac{1}{2}}{(1-p+b(\tau(\sigma))y^{2})}-p\frac{2\alpha}{1-p+\beta(\sigma) z^{2}},$$
the function
$$\mathcal{F}(a,b):=e^{\frac{ay^{2}}{4}}e^{-\frac{\alpha z^{2}}{4}}(\frac{\lambda_{1}}{\lambda})^{\frac{2p}{p-1}}F(a,b),$$
the nonlinear term
\begin{equation}\label{eq:DefNonlin}
\begin{array}{lll}
\mathcal{N}_{1}(a(\tau(\sigma)),b(\tau(\sigma)),\alpha,\eta(\sigma)):=(\frac{\lambda_{1}}{\lambda})^{\frac{2p}{p-1}}e^{\frac{a(\tau(\sigma))y^2}{4}}e^{-\frac{\alpha
z^{2}}{4}}N(a(\sigma),b(\tau(\sigma)),\xi(\tau(\sigma))).
\end{array}
\end{equation} where,
recall the definitions of $F$ and $N$ after (~\ref{eq:scource}),

In the next we prove that the new trajectory is a good approximation
of the old one.
\begin{proposition}\label{NewTrajectory}
For any $\tau\leq T$ we have that if $A(\tau)\leq
\beta^{-1/4}(\tau)$ then
\begin{equation}\label{eq:compare}
|\frac{\lambda}{\lambda_{1}}(t(\tau))-1|\lesssim \beta(\tau)
\end{equation}
for some constant $c$ independent of $\tau$.
\end{proposition}
\begin{proof} By the properties of $\lambda$ and $\lambda_{1}$ we have
\begin{equation}\label{EstLambda}
\frac{d}{d\tau}[\frac{\lambda}{\lambda_{1}}(t(\tau))-1]
=2a(\tau)(\frac{\lambda}{\lambda_{1}}(t(\tau))-1)+G(\tau)
\end{equation} with
$$G:=\alpha-a+(\alpha-a)(\frac{\lambda}{\lambda_{1}}-1)[(\frac{\lambda}{\lambda_{1}})^{2}+\frac{\lambda}{\lambda_{1}}+1]+a(\frac{\lambda}{\lambda_{1}}-1)^{2}[\frac{\lambda}{\lambda_{1}}+2].$$ By the definition of $A(\tau)$ we have that in the time interval
$\tau\in [0,T]$, if $A(\tau)\leq \beta^{-1/4}(\tau)$ then
\begin{equation}\label{CauchA}
|a(\tau)-\alpha|,\ |a(\tau)-\frac{1}{2}|\lesssim \beta(\tau).
\end{equation} Thus
\begin{equation}\label{Rem}
|G|\lesssim
\beta+(\frac{\lambda}{\lambda_{1}}-1)^{2}+|\frac{\lambda}{\lambda_{1}}-1|^{3}+\beta|\frac{\lambda}{\lambda_{1}}-1|.
\end{equation}
Observe that $\frac{\lambda}{\lambda_{1}}(t(\tau))-1=0$ when
$\tau=T.$ Thus Equations (~\ref{EstLambda}) can be rewritten as
\begin{equation}\label{LambdaRe}
\frac{\lambda_{1}}{\lambda}(t(\tau))-1=-\int_{\tau}^{T}e^{-\int^{s}_{\tau}2a(t)dt}G(s)ds.
\end{equation}
We claim that Equations (~\ref{CauchA}) and (~\ref{Rem}) imply
(~\ref{eq:compare}). Indeed, define an estimating function
$\Lambda(\tau)$ as
$$\Lambda(\tau):=\sup_{\tau\leq s\leq T}\beta^{-1}(s)|\frac{\lambda}{\lambda_{1}}(t(s))-1|.$$
Then (~\ref{LambdaRe}) and the assumption $A(\tau),\ B(\tau)\leq
\beta^{-1/4}(\tau)$ implies
$$
\begin{array}{lll}
|\frac{\lambda}{\lambda_{1}}(t(\tau))-1| &\lesssim &
\int_{\tau}^{T}e^{-\frac{1}{2}(T-\tau)}[\beta(s)+\beta^{2}(s)\Lambda^{2}(\tau)+\beta^{2}(s)\Lambda(\tau)]ds\\
&\lesssim
&\beta(\tau)+\beta^{2}(\tau)\Lambda^{2}(\tau)+\beta^{3}(\tau)\Lambda^{3}(\tau)+\beta^{2}(\tau)\Lambda(\tau),
\end{array}
$$
or equivalently
$$\beta^{-1}(\tau)|\frac{\lambda}{\lambda_{1}}(t(\tau))-1|\lesssim
1+\beta(\tau)\Lambda^{2}(\tau)+\beta^{2}(\tau)\Lambda^{3}(\tau)+\beta(\tau)\Lambda(\tau).$$
Consequently by the fact that $\beta(\tau)$ and $\Lambda(\tau)$ are
decreasing functions we have $$\Lambda(\tau)\lesssim
1+\beta(\tau)\Lambda^{2}(\tau)+\beta^{2}(\tau)\Lambda^{3}(\tau)+\beta(\tau)\Lambda(\tau)$$
which together with $\Lambda(T)=0$ implies $\Lambda(\tau)\lesssim 1$
for any time $\tau\in [0,T].$ Thus we prove the claim by the
definition of $\Lambda(\tau).$
\end{proof}
\begin{lemma}
For any $c_{1},c_{2}>0$ there exists a constant $c(c_{1},c_{2})$
such that
\begin{equation}\label{INT}
\int_{0}^{S}e^{-c_{1}(S-\sigma)}\beta^{c_{2}}(\tau(\sigma))d\sigma\leq
c(c_{1},c_{2})\beta^{c_{2}}(T).
\end{equation}
\end{lemma}
\begin{proof}
By the definition of $\tau(\sigma)$ we have that
$\sigma=\int_{0}^{t(\tau)}\lambda_{1}^{-2}(k)dk$ and
$\tau=\tau(\sigma)=\int_{0}^{t(\tau)}\lambda^{-2}(k)dk.$ By
Proposition ~\ref{NewTrajectory} we have that $\frac{1}{2}\leq
\frac{\lambda}{\lambda_{1}}\leq 2$, thus
\begin{equation}\label{TauS1}
\frac{1}{4}\sigma\geq \tau(\sigma)\geq 4\sigma
\end{equation} which
implies
$\frac{1}{\frac{1}{b_{0}}+\tau(\sigma)}\leq
\frac{4}{\frac{1}{b_{0}}+\sigma}.$ Moreover this yields
\begin{equation}\label{InI2}
\int_{0}^{S}e^{-c_{1}(S-s)}\beta^{c_{2}}(\tau(s))ds\leq
c(c_{1},c_{2})\frac{1}{(\frac{1}{b_{0}}+S)^{c_{2}}}.
\end{equation} Using
(~\ref{TauS1}) again to get $4S\geq \tau(S)=T\geq \frac{1}{4}S$
which together with (~\ref{InI2}) implies that
$$\int_{0}^{S}e^{-c_{1}(S-s)}\beta^{c_{2}}(\tau(s))ds\leq
c(c_{1},c_{2})\frac{1}{(\frac{1}{b_{0}}+T)^{c_{2}}}\lesssim
c(c_{1},c_{2})\beta^{c_{2}}(T).$$ Hence
$$\int_{0}^{S}e^{-c_{1}(S-s)}\beta^{c_{2}}(\tau(s))ds\leq
c(c_{1},c_{2})\beta^{c_{2}}(T)$$ which is (~\ref{INT}).
\end{proof}
\section{Propagator Estimates}\label{SEC:PropagatorEst}
In the next we present the decay estimates of the propagators
generated by $-L_{\alpha}$ and $-\mathcal{L}_{\alpha,\beta}$ which
play essential roles in proving the estimates for $M_{1}$, $M_{2}$
and $M_{q}$ below.

We start with analyzing the spectrum of the linearized operator
$\mathcal{L}_{\alpha,\beta}$ and $L_{\alpha}$. Due to the quadratic
term $\frac{1}{4} \alpha^{2} z^2$, the operators have discrete spectrum.
Then $L_{\alpha}+\frac{2p\alpha}{1-p}$ and $L_{\alpha}$ approximate
$\mathcal{L}_{\alpha\beta}$ near zero and at infinity respectively.
The spectrum of the operator $L_{\alpha}$ is
\begin{equation}\label{eq:spectrum}
\sigma(L_{\alpha})=\left\{n \alpha- \frac{2\alpha}{1-p}|\
n=0,1,2,\ldots\right\}.
\end{equation}
The first three normalized eigenvectors of $L_{\alpha}$, which are
used below, are
\begin{align}\label{eq:eigenvectors}
\phi_{0a}:=(\frac{\alpha}{2\pi})^\frac{1}{4} e^{-\frac{\alpha}{4}z^2},\
\phi_{1a}:=(\frac{\alpha}{2\pi})^{\frac{1}{4}}\sqrt{\alpha}z
e^{-\frac{\alpha}{4}z^2}, \phi_{2a}:=(\frac{\alpha}{8\pi})^{\frac{1}{4}}(1-\alpha
z^2)e^{-\frac{\alpha}{4}z^2}.
\end{align}
Denote the integral kernel of $e^{-\frac{\alpha
z^{2}}{4}}e^{-L_{\alpha}\sigma}e^{\frac{\alpha z^{2}}{4}}$ by
$U_{0}(x,y)$. By a standard formula (see \cite{Simon, GlJa}) we have
\begin{equation*}
e^{-L_{\alpha} \sigma}(x,y)=4\pi (1-e^{-2\alpha
\sigma})^{-1/2}\sqrt{\alpha}e^{\frac{2\alpha}{1-p}
\sigma}e^{-\alpha\frac{(x-e^{-\alpha \sigma}y)^{2}}{2(1-e^{-2\alpha
\sigma})}}. \label{eqn:96a}
\end{equation*}

In what follows we need the following result.
\begin{lemma}\label{kernelEst}
For any function $g$ and $\sigma>0$ we have that
\begin{equation}\label{est:99a}
\|\langle z\rangle^{-n}e^{\frac{\alpha
z^{2}}{4}}e^{-L_{\alpha}\sigma}g\|_{\infty}\lesssim e^{
\frac{2\alpha}{1-p} \sigma}\|\langle z\rangle^{-n}e^{\frac{\alpha
z^{2}}{4}}g\|_{\infty},\ n=0,1,2,
\end{equation}
 or equivalently
\begin{equation}\label{eq:secondForm}
e^{\frac{\alpha x^{2}}{2}}\int \langle
x\rangle^{-2}U_{0}(x,y)e^{-\frac{\alpha}{2}y^{2}}\langle
y\rangle^{2}dy\lesssim e^{\frac{2\alpha}{1-p} \sigma}.
\end{equation}
\end{lemma}
\begin{proof}
Note that the first three eigenvectors of $L_{\alpha}$ are
$e^{-\frac{\alpha z^{2}}{4}},\ z e^{-\frac{\alpha z^{2}}{4}},\
(\alpha z^{2}-1)e^{-\frac{\alpha z^{2}}{4}}$ with the eigenvalues
$-\frac{2\alpha}{1-p}, -\frac{2\alpha}{1-p}+\alpha,\
-\frac{2p\alpha}{1-p}$ (see (~\ref{eq:spectrum})). Using that the
integral kernel of $e^{-\sigma L_{\alpha}}$ is positive and
therefore $\|e^{-\sigma L_{\alpha}} g\|_\infty\le \|f^{-1}
g\|_\infty\|e^{-\sigma L_{\alpha}} f\|_\infty$ for any $f>0$ and
using that $e^{-\sigma L_{\alpha}} e^{-\frac{\alpha}{4}
z^2}=e^{\frac{2\alpha}{1-p} \sigma}e^{-\frac{\alpha}{4}z^2}$ and
$e^{-\sigma L_{\alpha}}(\alpha z^2-1)e^{-\frac{\alpha}{4}
z^2}=e^{\frac{2p\alpha}{1-p} \sigma}(\alpha
z^2-1)e^{-\frac{\alpha}{4}z^2}$, we find that
$$
\begin{array}{lll}
\|\langle z\rangle^{-2}e^{\frac{\alpha z^{2}}{4}}e^{-\sigma
L_{\alpha}}g\|_{\infty} &\leq &\|\langle z\rangle^{-2}
e^{\frac{\alpha z^{2}}{4}}e^{-\sigma L_{\alpha}}e^{-\frac{\alpha
z^{2}}{4}}(z^{2}+1)\|_{\infty}\|\langle z\rangle^{-2}e^{\frac{\alpha
z^{2}}{4}}g\|_{\infty}\\
&=&\|\langle z\rangle^{-2}[e^{\frac{2\alpha}{1-p}
\sigma}\frac{1}{\alpha}+e^{\frac{2p\alpha}{1-p}\sigma}(z^{2}-\frac{1}{\alpha})]\|_{\infty}\|\langle
z\rangle^{-2}e^{\frac{\alpha
z^{2}}{4}}g\|_{\infty}\\
&\leq&2(\frac{1}{\alpha}+1)e^{\frac{2\alpha}{1-p} \sigma}\|\langle
z\rangle^{-2}e^{\frac{\alpha z^{2}}{4}}g\|_{\infty}
\end{array}
$$ and $$
\begin{array}{lll}
\|e^{\frac{\alpha z^{2}}{4}}e^{-\sigma L_{\alpha}}g\|_{\infty} &\leq
&\| e^{\frac{\alpha z^{2}}{4}}e^{-\sigma L_{\alpha}}e^{-\frac{\alpha
z^{2}}{4}}\|_{\infty}\|e^{\frac{\alpha
z^{2}}{4}}g\|_{\infty}\\
&=& e^{\frac{2\alpha}{1-p}}\|e^{\frac{\alpha z^{2}}{4}}g\|_{\infty}
\end{array}
$$ which are the case $n=0,2$ of (~\ref{est:99a}). The case $n=1$
follows from the interpolation between $n=0$ and $n=2.$ Moreover
this together with the definition of $U_{0}(x,y)$ after
(~\ref{eq:eigenvectors}) and setting $g=\langle z\rangle^{2}
e^{-\frac{\alpha z^{2}}{4}}$ implies \eqref{eq:secondForm}.
\end{proof}
We define $P^{\alpha}_{n}=1-\overline{P^{\alpha}_{n}}$ and define $\overline{P^{\alpha}_{n}},\ n=1,2,3,$ as the projection onto
the space spanned by the first $n$ eigenvectors of $L_{\alpha}$ with
the form
\begin{equation}\label{eq:projection}
\begin{array}{lll}
\overline{P^{\alpha}_{n}}&=&\displaystyle\sum_{m=0}^{n-1}| \phi_{m,\alpha}
\rangle \langle \phi_{m,\alpha} |,
\end{array}
\end{equation} where, recall the
definitions of $\phi_{m,\alpha}$ in (~\ref{eq:eigenvectors}).
\begin{proposition}\label{PRO:propagator}
Let $P^{\alpha}_{ n}$ be the projection defined above. Then for any function $g$ and time
$ \sigma\geq 0$ we have
\begin{equation}\label{eq:estproject2}
\|\langle z\rangle^{-2}e^{\frac{\alpha
z^{2}}{4}}e^{-L_{\alpha}\sigma}P_{2}^{\alpha}g\|_{\infty}\lesssim
e^{\frac{2\alpha p}{1-p} \sigma}\|\langle
z\rangle^{-2}e^{\frac{\alpha z^{2}}{4}}g\|_{\infty},
\end{equation}
\begin{equation}\label{eq:hPropagator}
\|\langle z\rangle^{-k}e^{\frac{\alpha
z^{2}}{4}}e^{-L_{\alpha}\sigma}P_{1}^{\alpha}g\|_{\infty}\lesssim
e^{[\frac{2}{1-p}-k]\alpha\sigma}\|\langle
z\rangle^{-k}e^{\frac{\alpha z^{2}}{4}}g\|_{\infty}
\end{equation} for any $k\in [0,1]$;
and there
exists constant $c_{0}>0$ such that for any $\tau\geq \sigma\geq 0$
\begin{equation}\label{eq:estProject3}
\|\langle z\rangle^{-3}e^{\frac{\alpha
z^{2}}{4}}P_{3}^{\alpha}U(\tau,\sigma)P_{3}^{\alpha}g\|_{\infty}\lesssim
e^{-c_{0}(\tau-\sigma)}\|\langle z\rangle^{-3}e^{\frac{\alpha
z^{2}}{4}}g\|_{\infty}
\end{equation} where $U(\tau,\sigma)$ denotes the propagator
generated by the operator
$-P_{3}^{\alpha}\mathcal{L}_{\alpha,\beta}P_{3}^{\alpha}$.
\end{proposition}
\begin{proof}
(~\ref{eq:estProject3}) is proved in ~\cite{DGSW}.

Now we prove (~\ref{eq:estproject2}). Define a new function
$f:=e^{-\frac{\alpha y^{2}}{4}}P^{\alpha}g$. The definition of
$U_{0}(x,y)$ after (~\ref{eq:eigenvectors}) implies
\begin{equation}\label{FK1}
e^{-L_{\alpha}\sigma}P^{\alpha}g=\int e^{\frac{\alpha
x^{2}}{4}}U_{0}(x,y)f(y)dy.
\end{equation}
Integrate by parts on the right hand side of (~\ref{FK1}) to obtain
\begin{equation}\label{Estimateta}
\begin{array}{lll}
e^{-L_{\alpha} \sigma}P^{\alpha}g=e^{\frac{\alpha x^{2}}{4}}\int
\partial_{y}^{2}U_{0}(x,y)f^{(-2)}(y)dy
\end{array}
\end{equation} where
$f^{(-m-1)}(x):=\int_{-\infty}^{x}f^{(-m)}(y)dy$ and $f^{(0)}:=f.$
Now we estimate the right hand side of Equation (~\ref{Estimateta}).
\begin{enumerate}
\item[(A)] By the facts that $f=e^{-\frac{\alpha
y^{2}}{4}}P^{\alpha}g$ and $P_{2}^{\alpha}g\perp
y^{n}e^{-\frac{\alpha y^{2}}{4}},\ n=0,1,$ we have that $f\perp
1,\ y.$ Therefore by integration by parts we have
$$f^{(-m)}(y)=\int_{-\infty}^{y}f^{(-m+1)}(x)dx=-\int_{y}^{\infty}f^{(-m+1)}(x)dx,\ m=1,2$$
which together with the definition of $f^{(-m)}$ yields
$$|f^{(-2)}(y)|\lesssim e^{-\frac{\alpha}{2}y^{2}}\|\langle
y\rangle^{-2}e^{\frac{\alpha}{4}y^{2}}P^{\alpha}g\|_{\infty}.$$
\item[(B)] Using the explicit formula for $U_0(x,y)$ given above we find
$$|\partial^{2}_{y}U_{0}(x,y)|\lesssim \frac{e^{-2\alpha
\sigma}}{(1-e^{-2\alpha \sigma})^{2}}(|x|+|y|)^{2}U_{0}(x,y).$$
\end{enumerate}

Collecting the estimates (A)-(B) above and using Equation
(~\ref{Estimateta}), we have
$$
\begin{array}{lll}
& &\langle x\rangle^{-2}e^{\frac{\alpha x^{2}}{4}}| e^{-L_{\alpha} \sigma}P^{\alpha}g(x)|\\
&\lesssim & \frac{e^{-2\alpha \sigma}}{(1-e^{-2\alpha
\sigma})^{2}}\langle
x\rangle^{-2}e^{\frac{\alpha x^{2}}{2}}\int (|x|+|y|)^{2}U_{0}(x,y)|f^{(-2)}(y)|dy\\
&\lesssim &\frac{e^{-2\alpha \sigma}}{(1-e^{-2\alpha
\sigma})^{2}}e^{\frac{\alpha x^{2}}{2}}\int \langle
x\rangle^{-2}U_{0}(x,y)e^{-\frac{\alpha}{2}y^{2}}\langle
y\rangle^{2}dy \|\langle
y\rangle^{-2}e^{\frac{\alpha}{4}y^{2}}P_{2}^{\alpha}g\|_{\infty}\\
&\lesssim &\frac{e^{-2\alpha \sigma}}{(1-e^{-2\alpha
\sigma})^{2}}e^{\frac{\alpha x^{2}}{2}}\int \langle
x\rangle^{-2}U_{0}(x,y)e^{-\frac{\alpha}{2}y^{2}}\langle
y\rangle^{2}dy \|\langle
y\rangle^{-2}e^{\frac{\alpha}{4}y^{2}}g\|_{\infty}
\end{array}
$$ where in the last step we used the explicit form of
$P_{2}^{\alpha}$ in (~\ref{eq:projection}). This together with the
estimate \eqref{eq:secondForm} gives the estimate
(~\ref{eq:estproject2}) when $\sigma>1$. If $\sigma\leq 1$ we use
(~\ref{est:99a}) to remove the singularity at $\sigma=0.$

When $k=1$ the proof of (~\ref{eq:hPropagator}) is almost the same
to the proof of (~\ref{eq:estproject2}) and, thus omitted; when $k=0$ we
have $$\|e^{\frac{\alpha
z^{2}}{4}}e^{-L_{\alpha}\sigma}P_{1}^{\alpha}g\|_{\infty}\lesssim
e^{\frac{2}{1-p}\alpha\sigma}\|e^{\frac{\alpha
z^{2}}{4}}P_{1}^{\alpha}g\|_{\infty}\lesssim \|e^{\frac{\alpha
z^{2}}{4}}g\|_{\infty}$$ by using (~\ref{est:99a}) and the observation
that $\|P_{1}g\|_{\infty}\lesssim \|g\|_{\infty}$. The general case
follows from the interpolation between $k=1$ and $k=0.$

Thus the proof is complete.
\end{proof}
\section{Estimate of $M_{1}$}\label{SEC:EstM1}
% Proofread\\
In this subsection we derive an estimate for $M_{1}$ in Equations
(~\ref{eq:majorants}).

Given any time $\tau$, choose $T=\tau$ and do the estimates as in
Proposition ~\ref{NewTrajectory}. We start from estimating $\eta$ in
Equation (~\ref{eq:eta}). Observe that the function $\eta$ is not
orthogonal to the first three eigenvectors of the operator
$L_{\alpha}$, defined in (~\ref{eq:eta}), thus we put projections on
both sides of Equation (~\ref{eq:eta}) to get
\begin{equation}\label{EQ:eta2}
\frac{d}{d\sigma}P_{3}^{\alpha}\eta=-P_{3}^{\alpha}\mathcal{L}_{\alpha\beta}P_{3}^{\alpha}\eta+\sum_{n=1}^{4}P_{3}^{\alpha}D_{n}
\end{equation}
where the functions $D_{n}, \ n=1,2,3,4,$ are defined as
$$D_{1}:=-P_{3}^{\alpha}V\eta+P_{3}^{\alpha}VP_{3}^{\alpha}\eta,\ \ \ D_{2}:=\mathcal{W}\eta,$$
$$D_{3}:=\mathcal{F}(a,b),\ \ \  D_{4}:=\mathcal{N}_{1}(a,b,\alpha,\eta),$$
where, recall the definitions of the functions $\mathcal{F},\
\mathcal{N}_{1}$ and the operators $\mathcal{W},\ V$ after Equation
(~\ref{eq:eta}).

Now we start with estimating the terms $D_{n}$, $n=1,2,3,4$, on the
right hand side of (~\ref{EQ:eta2}).
\begin{lemma}\label{EstDs} Suppose that $M_{1}(\tau)\leq \frac{1}{8}(1-p)^{\frac{1}{1-p}-\frac{3}{2}},\ A(\tau),\ B(\tau)\leq
\beta^{-1/4}(\tau)$ and the function $v(y,\tau)$ defined in
(~\ref{eq:definev}) satisfies the estimate $v(y,\tau)\geq
\frac{1}{2}(\frac{1-p}{2})^{\frac{1}{1-p}}$. Then we have
\begin{equation}\label{eq:estD1}
\|\langle z\rangle^{-3}e^{\frac{\alpha
z^{2}}{4}}D_{1}(\sigma)\|_{\infty}\lesssim
\beta^{2}(\tau(\sigma))M_{1}(T),
\end{equation}
\begin{equation}\label{est:termW3}
\|\langle z\rangle^{-3}e^{\frac{\alpha
z^{2}}{4}}D_{2}(\sigma)\|_{\infty}\lesssim
\beta^{7/4}(\tau(\sigma))M_{1}(T),
\end{equation}
\begin{equation}\label{eq:estSource3}
\|\langle z\rangle^{-3}e^{\frac{\alpha
z^{2}}{4}}D_{3}(\sigma)\|_{\infty} \lesssim
 \beta^{5/2}(\tau(\sigma))[1+M^{2-p}_{1}(T)+A(T)M_{1}(T)],
\end{equation}
if $ -1\leq p<0$ then
\begin{equation}\label{eq:nonlinearity3}
\|\langle z\rangle^{-3}e^{\frac{\alpha
z^{2}}{4}}D_{4}(\sigma)\|_{\infty}\lesssim
\beta^{3/2}(\tau(\sigma))M_{1}(T)M_{2}(T).
\end{equation}
if $p<-1$ then
\begin{equation}\label{eq:hnonlinearity3}
\|\langle z\rangle^{-3}e^{\frac{\alpha
z^{2}}{4}}D_{4}(\sigma)\|_{\infty}\lesssim
\beta^{3/2}(\tau(\sigma))M_{1}(T)[M_{q}^{1-p}(T)+M_{q}(T)]
\end{equation} where, recall that $q=\min\{\frac{4}{1-p},\frac{2(2-p)}{(p-1)^2},1\}.$
\end{lemma}
\begin{proof}
In what follows we implicitly use
\begin{equation}\label{eq:compare2}
\frac{\lambda_{1}}{\lambda}(t(\tau))-1=O(\beta(\tau)), \
\text{thus}\ \frac{\lambda_{1}}{\lambda}(t(\tau)),\
\frac{\lambda}{\lambda_{1}}(t(\tau))\leq 2
\end{equation} implied by (~\ref{eq:compare}) and the assumptions on $M_{1}, A,\ B$.

By the relation between $\xi$ and $\eta$ in (~\ref{NewFun}) and the
fact that $y=\frac{\lambda_{1}}{\lambda}z$, we have
$$\|\langle z\rangle^{-n}e^{\frac{\alpha
z^{2}}{4}}\eta(\sigma)\|_{\infty}\lesssim \|\langle
y\rangle^{-n}e^{\frac{a y^{2}}{4}}\xi(\tau(\sigma))\|_{\infty}, $$
$n=2,3$ and $n=q$  for $p<-1.$ By the definition of $M_{1},\ M_{2}$
and the fact that $\tau(\sigma)\leq T$ we have
\begin{equation}\label{eq:compareXiEtaN} \|\langle
z\rangle^{-n}e^{\frac{\alpha
z^{2}}{4}}\eta(\sigma)\|_{\infty}\lesssim
\beta^{\frac{n}{2}}(\tau(\sigma))M_{k_{n}}(T)
\end{equation} with $n=2,3$, especially $n=q$ for $p<-1,$ and $k_{2}=2,\ k_{3}=3,\ k_{q}=q$. This proves (~\ref{eq:estD1}) as a special case.

We rewrite $D_{1}=P_{3}^{\alpha}D_{1}$ as
$$P_{3}^{\alpha}D_{1}(\sigma)=-P_{3}^{\alpha}p\frac{\alpha+\frac{1}{2}}{p-1+b(\tau(\sigma))z^{2}}b(\tau(\sigma))z^{2}(1-P_{3}^{\alpha})\eta(\sigma)$$
which admits the estimate $$
\begin{array}{lll}
\|\langle z\rangle^{-3}e^{\frac{\alpha
z^{2}}{4}}P_{3}^{\alpha}D_{1}(\sigma)\|_{\infty}&\lesssim& |\langle
z\rangle^{-1}\frac{b(\tau(\sigma))z^{2}}{1+bz^{2}}|\| \langle
z\rangle^{-2}e^{\frac{\alpha z^{2}}{4}}(1-P_{3}^{\alpha})\eta(\sigma)\|_{\infty}\\
&\lesssim & b^{1/2}(\tau(\sigma))\|\langle
z\rangle^{-3}e^{\frac{\alpha
z^{2}}{4}}\eta(\sigma)\|_{\infty}\\
&\lesssim & \beta^{1/2}(\tau(\sigma))\|\langle
z\rangle^{-3}e^{\frac{\alpha z^{2}}{4}}\eta(\sigma)\|_{\infty}
\end{array}
$$ where we use that $b(\tau)\lesssim \beta(\tau)$ implied by $B(\tau)\leq
\beta^{-1/4}(\tau)$ and the fact that
\begin{equation}\label{eq:estPro}
\|\langle z\rangle^{-2}e^{\frac{\alpha
z^{2}}{4}}(1-P_{3}^{\alpha})g\|_{\infty}\lesssim \|\langle
z\rangle^{-3}e^{\frac{\alpha z^{2}}{4}}g\|_{\infty}
\end{equation} from the explicit form of $1-P_{3}^{\alpha}$ in (~\ref{eq:projection}). This estimate together with (~\ref{eq:compareXiEtaN}) implies (~\ref{eq:estD1}).

Now we prove (~\ref{est:termW3}). Recall that
$y=\frac{\lambda_{1}}{\lambda}z.$ After some manipulation on the expression of $\mathcal{W}$ we have
$$|D_{2}(\sigma)|\lesssim
(|\frac{\lambda}{\lambda_{1}}-1|+|a(\tau(\sigma))-\alpha|+|\frac{b(\tau(\sigma))-\beta(\tau(\sigma))}{\beta(\tau(\sigma))}|)|\eta(\sigma)|.$$
Equations (~\ref{eq:compare}) and (~\ref{CauchA}) imply that
$|\frac{\lambda}{\lambda_{1}}-1|+|a(\tau(\sigma))-\alpha|\lesssim
\beta(\tau(\sigma))$; the assumption on $B$ and its definition imply
$|\frac{b(\tau(\sigma))-\beta(\tau(\sigma))}{\beta(\tau(\sigma))}|\lesssim
\beta^{1/4}(\tau(\sigma)).$ Consequently
$$\|\langle z\rangle^{-3}e^{\frac{\alpha z^{2}}{4}}D_{2}(\sigma)\|\lesssim \beta^{1/4}(\tau(\sigma))\|\langle z\rangle^{-3}e^{\frac{\alpha
z^{2}}{4}}\eta(\sigma)\|_{\infty}$$ which together with
(~\ref{eq:compareXiEtaN}) implies Equation (~\ref{eq:estD1}).

For (~\ref{eq:estSource3}), by the relation between
$D_{3}$, $\mathcal{F}$, $F$ in (~\ref{eq:eta}), (~\ref{EQ:eta2}) and
Equation (~\ref{eq:compare2}) we have
$$\|\langle z\rangle^{-3}e^{\frac{\alpha z^{2}}{4}}D_{3}(\sigma)\|_{\infty}\lesssim \|\langle
y\rangle^{-3}e^{\frac{a(\tau(\sigma))y^{2}}{4}}F(a(\tau(\sigma)),b(\tau(\sigma)))\|_{\infty}$$
which together with the estimate of $F(a,b)$ in
(~\ref{eq:estSource}) implies (~\ref{eq:estSource3}).

When $0>p\geq -1$, by the relation between $D_{4}$,
$\mathcal{N}_{1}(a,b,\alpha,\eta)$ and $N(a,b,\xi)$ in
(~\ref{EQ:eta2}), (~\ref{eq:DefNonlin}) and the estimate
(~\ref{eq:ESTnonlin}) we have
$$|e^{\frac{\alpha z^{2}}{4}}D_{4}(\sigma)|\lesssim (\frac{1}{1+\beta(\tau(\sigma)) y^{2}})^{\frac{2}{1-p}}e^{\frac{ay^{2}}{2}}|\xi|^{2}\leq \frac{1}{\beta(\tau(\sigma))}\langle y\rangle^{-2}e^{\frac{ay^{2}}{2}}|\xi|^{2},$$
which together with (~\ref{eq:compareXiEtaN}) implies
(~\ref{eq:nonlinearity3}).

When $p<-1$, by (~\ref{eq:pESTnonlin}) and the definition of $q$ we
have
$$
\begin{array}{lll}
|e^{\frac{\alpha z^{2}}{4}}D_{4}(\sigma)|&\lesssim
&(\frac{1}{1+\beta
y^{2}})^{\frac{q(1-p)}{2}}|e^{\frac{ay^{2}}{4}}\xi|^{2-p}+(\frac{1}{1+\beta
y^{2}})^{\frac{q}{2}}|e^{\frac{ay^{2}}{4}}\xi|^{2}\\
&\lesssim& \beta^{\frac{q(p-1)}{2}}\langle
y\rangle^{q(p-1)}|e^{\frac{ay^{2}}{4}}\xi|^{2-p}+\beta^{-\frac{q}{2}}\langle
y\rangle^{-q}|e^{\frac{ay^{2}}{4}}\xi|^{2}
\end{array}
$$ which together with
(~\ref{eq:compareXiEtaN}) implies (~\ref{eq:hnonlinearity3}).

Hence the proof is complete.
\end{proof}
\subsection{Proof of Equations (~\ref{eq:M1}) and (~\ref{eq:hM1})}
% Proofread\\
By Duhamel principle we rewrite Equation (~\ref{EQ:eta2}) as
\begin{equation}\label{eq:eta3}
P_{3}^{\alpha}\eta(S)=U(S,0)P_{3}^{\alpha}\eta(0)+\displaystyle\sum_{n=1}^{4}\int_{0}^{S}U(S,\sigma)P_{3}^{\alpha}D_{n}(\sigma)d\sigma,
\end{equation} where, recall, $U(t,s)$ is the propagator
generated by the operator
$-P_{3}^{\alpha}\mathcal{L}_{\alpha,\beta}P_{3}^{\alpha}$. Use
(~\ref{eq:estProject3}) to get
\begin{equation}\label{eq:M1Ge}
\begin{array}{lll}
& &\beta^{-3/2}(T)\|\langle z\rangle^{-3}e^{\frac{\alpha
z^{2}}{4}}P_{3}^{\alpha}\eta(S)\|_{\infty} \lesssim X_{1}+X_{2}
\end{array}
\end{equation} with $$X_{1}:=e^{-c_{0}S}\beta^{-3/2}(T)\|\langle
z\rangle^{-3}e^{\frac{\alpha z^{2}}{4}}\eta(0)\|_{\infty},$$
$$X_{2}:=\beta^{-3/2}(T)\displaystyle\int_{0}^{S}e^{-c_{0}(S-\sigma)}\sum_{n=1}^{4}\|\langle
z\rangle^{-3}e^{\frac{\alpha
z^{2}}{4}}D_{n}(\sigma)\|_{\infty}d\sigma.$$

Now we estimate each term on the right hand side.
\begin{enumerate}
\item[(A)] The slow
decay of $\beta(\tau)$ implies $e^{-c_{0}S}\beta^{-3/2}(T)\lesssim \beta^{-3/2}(0)$. Then we use Equation (~\ref{eq:compareXiEtaN}) to obtain
\begin{equation}
\begin{array}{lll}
X_{1}\lesssim \beta^{-3/2}(0)\|\langle z\rangle^{-3}e^{\frac{\alpha
z^{2}}{4}}\eta(0)\|_{\infty} \lesssim M_{1}(0).
\end{array}
\end{equation}
\item[(B)] By the integral estimate (~\ref{INT}), the estimates of $D_{n},\ n=1,2,3,4,$ in
(~\ref{eq:estD1})-(~\ref{eq:nonlinearity3}) and the fact
$\beta(\tau)\leq \beta(0)$ we obtain
\begin{equation}\label{eq:estD3}
X_{2} \lesssim
\beta^{-3/2}(T)\int_{0}^{S}e^{-c_{0}(S-\sigma)}\beta^{3/2}(\tau(\sigma))
d\sigma X^{(p)}\lesssim X^{(p)}
\end{equation} with
$$X^{(p)}:=
\left\{
\begin{array}{lll}
\beta^{1/4}(0)[1+M_{1}^{2-p}(T)+A(T)M_{1}(T)]+M_{1}(T)M_{2}(T),\
& \text{if}\ 0>p\geq -1;\\
\beta^{1/4}(0)[1+M_{1}^{2-p}(T)+A(T)M_{1}(T)]+M_{1}(T)[M_{q}^{1-p}(T)+M_{q}(T)],\
& \text{if}\ p<-1.
\end{array}
\right.
$$
\end{enumerate}
In Equation (~\ref{NewFun}) we define
$\lambda_{1}(t(T))=\lambda(t(T))$, $\xi(\cdot, T)=\eta(\cdot,S)$ and
$\alpha=a(T)$, thus $z=y$ and $P^{\alpha}_{3}\eta(S)=\xi(T)$, consequently
$$\|\langle z\rangle^{-3}e^{\frac{\alpha
z^{2}}{4}}P_{3}^{\alpha}\eta(S)\|_{\infty}=\|\langle
y\rangle^{-3}e^{\frac{ay^{2}}{4}}\xi(T)\|_{\infty}$$ which together
with Equations (~\ref{eq:M1Ge})-(~\ref{eq:estD3}) implies
$$
\beta^{-3/2}(T)\|\langle
y\rangle^{-3}e^{\frac{ay^{2}}{4}}\xi(T)\|_{\infty} \lesssim
M_{1}(0)+X^{(p)}+\beta^{1/4}(0)[1+M_{1}^{2-p}(T)+A(T)M_{1}(T)].
$$

By the definition of $M_{1}$ in (~\ref{eq:majorants}), we obtain
$$
M_{1}(T)\lesssim
M_{1}(0)+X^{(p)}+\beta^{1/4}(0)[1+M_{1}^{2-p}(T)+A(T)M_{1}(T)].
$$
This together with the fact that $T$ is arbitrary implies Equation
(~\ref{eq:M1}) and (~\ref{eq:hM1}).
\begin{flushright}
$\square$
\end{flushright}

%%%%%%%%%%%%%%%%%%%%%%%%%%%%%%%%%%%%%%%%%%%%%%%%%%%%%%%%%%%%%%%%%%%%%%%%%%%%%%%%
%%%%%%%%%%%%%%%%%%%%%%%%%%%%%%%%%%%%%%%%%%%%%%%%%%%%%%%%%%%%%%%%%%%%%%%%%%%%
\section{Proof of Equations (~\ref{eq:M2}) and (~\ref{eq:hM2})}\label{SEC:EstM2}
We rewrite Equation (~\ref{eq:eta}) as
\begin{equation}\label{eq:eta4}
\frac{d}{d\sigma}\eta(\sigma)=-L_{\alpha}\eta(\sigma)-V\eta+\sum_{n=2}^{4}D_{n}
\end{equation}
where, recall the definitions of the operators $L_{\alpha},\ -V$ in
(~\ref{eq:eta}), the definition of $D_{n},\ n=2,3,4,$ in
(~\ref{EQ:eta2}).
\begin{lemma}
Suppose that $M_{1}(\tau)\leq \frac{1}{8}(1-p)^{\frac{1}{1-p}-\frac{3}{2}},\
A(\tau),\ B(\tau)\leq \beta^{-1/4}(\tau)$ and the function
$v(y,\tau)$ defined in (~\ref{eq:definev}) has the lower bound
$v(y,\tau)\geq \frac{1}{2}(\frac{1-p}{2})^{\frac{1}{1-p}}$. Then we
have
\begin{equation}\label{eq:estF3}
\|\langle z\rangle^{-2}e^{\frac{\alpha
z^{2}}{4}}V\eta(\sigma)\|_{\infty}\lesssim \beta(\tau(\sigma))
M_{1}(T),
\end{equation}
\begin{equation}\label{est:termW2}
\|\langle z\rangle^{-2}e^{\frac{\alpha
z^{2}}{4}}D_{2}(\sigma)\|_{\infty}\lesssim
\beta^{5/4}(\tau(\sigma))M_{2}(T),
\end{equation}
\begin{equation}\label{eq:estSource2}
\|\langle z\rangle^{-2}e^{\frac{\alpha
z^{2}}{4}}D_{3}(\sigma)\|_{\infty} \lesssim
\beta^{2}(\tau(\sigma))[1+M^{2-p}_{1}(T)+A(T)M_{1}(T)],
\end{equation}
if $0<p\leq -1$ then
\begin{equation}\label{eq:pnonlinearity2}
\|\langle z\rangle^{-2}e^{\frac{\alpha
z^{2}}{4}}D_{4}(\sigma)\|_{\infty}\lesssim
\beta(\tau(\sigma))M_{2}^{2}(T);
\end{equation}
if $p<-1$ then
\begin{equation}\label{eq:nonlinearity2}
\|\langle z\rangle^{-2}e^{\frac{\alpha
z^{2}}{4}}D_{4}(\sigma)\|_{\infty}\lesssim
\beta(\tau(\sigma))M_{2}(T)[M_{q}^{1-p}(T)+M_{q}(T)].
\end{equation}
\end{lemma}
\begin{proof}
The proof is easier than that of Lemma ~\ref{EstDs}, thus is
omitted.
\end{proof}
Rewrite (~\ref{eq:eta4}) to have
$$
P_{2}^{\alpha}\eta(S)=e^{-L_{\alpha}S}P_{2}^{\alpha}\eta(0)+\int_{0}^{S}e^{-L_{\alpha}(S-\sigma)}P_{2}^{\alpha}[-V\eta(\sigma)+\sum_{n=2}^{4}D_{n}]d\sigma,
$$ where, recall the definition
of $S$ in (~\ref{T2}). The propagator estimate of
$e^{-L_{\alpha}\sigma}P_{2}^{\alpha}$ in (~\ref{eq:estproject2})
yields
\begin{equation}\label{K123s}
\begin{array}{lll}
\beta^{-1}(T)\|\langle z\rangle^{-2}e^{\frac{\alpha
z^{2}}{4}}P_{2}^{\alpha}\eta(S)\|_{\infty}\lesssim K_{0}+K_{1}
\end{array}
\end{equation} where $K_{n}$'s are given by
$$K_{0}:=e^{\frac{2p\alpha}{1-p} S}\beta^{-1}(T)\|\langle
z\rangle^{-2}e^{\frac{\alpha z^{2}}{4}}\eta(0)\|_{\infty},$$
$$K_{1}:=\beta^{-1}(T)\int_{0}^{S}e^{\frac{2p\alpha}{1-p}(S-\sigma)}[\|\langle z\rangle^{-2}e^{\frac{\alpha z^{2}}{4}}V\eta(\sigma)\|_{\infty}+\sum_{n=2}^{4}\|\langle z\rangle^{-2}e^{\frac{\alpha
z^{2}}{4}}D_{n}(\sigma)\|_{\infty}]d\sigma,$$ where, recall that
$p<0.$

In the next we estimate $K_{n}$'s, $n=0,1.$
\begin{itemize}
\item[(K0)] First, $K_{0}$ has the bound
$$K_{0}\lesssim \beta^{-1}(T)e^{\frac{2p\alpha}{1-p} S}\|\langle z\rangle^{-2}e^{\frac{\alpha z^{2}}{4}}\eta(0)\|_{\infty}$$
The slow decay of $\beta$ and Equation (~\ref{eq:compareXiEtaN})
yield
\begin{equation}\label{eq:estK0}
K_{0}\lesssim \beta^{-1}(0)\|\langle z\rangle^{-2}e^{\frac{\alpha
z^{2}}{4}}\eta(0)\|_{\infty}\lesssim M_{2}(0).
\end{equation}
\item[(K1)]
By the estimates of $V\eta,\ D_{2},\ D_{3},\ D_{4},$ in Equations
(~\ref{eq:estF3})-(~\ref{eq:nonlinearity2}) and the integral
estimate (~\ref{INT}), we have that if $0>p\geq -1$ then
\begin{equation}\label{eq:estK2}
K_{2}\lesssim K^{(p)}
\end{equation} with the constant
$$K^{(p)}:=
M_{1}(T)+M_{2}^{2}(T)+\beta^{1/4}(0)[1+M_{1}^{2-p}(T)+A(T)M_{1}(T)]\
\text{if}\ 0>p\geq -1;$$ and $$
K^{(p)}:=M_{1}(T)+M_{2}(T)(M_{q}^{1-p}(T)+M_{q}(T))+\beta^{1/4}(0)[1+M_{1}^{2-p}(T)+M_{2}(T)+A(T)M_{1}(T)]$$
if $p<-1 . $
\end{itemize}
Collecting the estimates (~\ref{K123s})-(~\ref{eq:estK2}) we have
\begin{equation}\label{FinalStep}
\begin{array}{lll}
\beta^{-3/2}(T)\|\langle z\rangle^{-2}e^{\frac{\alpha
z^{2}}{4}}P_{2}^{\alpha}\eta(S)\|_{\infty} \lesssim
M_{2}(0)+K^{(p)}.
\end{array}
\end{equation}
Recall the definition of $T$ in (~\ref{NewFun}). By the relation
between $\xi$ and $\eta$ in Equation (~\ref{NewFun}), we have
$\xi(T)=\eta(S),$ $\alpha=a(T)$, and $y=z$, hence
$P_{2}^{\alpha}\eta(S)=\xi(T)$ and
$$
\begin{array}{lll}
\beta^{-3/2}(T)\|\langle
y\rangle^{-2}e^{\frac{ay^{2}}{4}}\xi(T)\|_{\infty}
&=&\beta^{-3/2}(T)\|\langle z\rangle^{-2}e^{\frac{\alpha
z^{2}}{4}}P_{2}^{\alpha}\eta(S)\|_{\infty}
\end{array}
$$
which together with (~\ref{FinalStep}) implies
$$
\begin{array}{lll}
M_{2}(T)&\lesssim & M_{2}(0)+K^{(p)}.
\end{array}
$$
Since $T$ is an arbitrary time, the proofs of (~\ref{eq:M2}) and
(~\ref{eq:hM2}) are complete.

%%%%%%%%%%%%%%%%%%%%%%%%%%%%%%%%%%%%%%%%%%%%%%%%%%%%%%%%%%%%%%%%%%%%%%%%%%%%%%%%
%%%%%%%%%%%%%%%%%%%%%%%%%%%%%%%%%%%%%%%%%%%%%%%%%%%%%%%%%%%%%%%%%%%%%%%%%%%%%
\section{Proof of Equation (~\ref{eq:M3})}\label{SEC:EstM3}
% Proofread\\
In the following lemma we present the estimates for $D_{n},\
n=2,3,4$ in (~\ref{eq:eta4}). Recall that $q=\min\{\frac{4}{1-p},\
\frac{2(2-p)}{(1-p)^{2}},\ 1\}.$
\begin{lemma}
Suppose that for $\tau\leq T$, $M_{1}(\tau)\leq
\frac{1}{8}(1-p)^{\frac{1}{1-p}-\frac{3}{2}},\ A(\tau),\ B(\tau)\leq
\beta^{-1/4}(\tau)$ and the function $v(y,\tau)$ defined in
(~\ref{eq:definev}) has the lower bound $v(y,\tau)\geq
\frac{1}{2}(\frac{1-p}{2})^{\frac{1}{1-p}}$. Then we have
\begin{equation}\label{eq:hestF3}
\|\langle z\rangle^{-q}e^{\frac{\alpha
z^{2}}{4}}V\eta(\sigma)\|_{\infty}\lesssim
\beta^{\frac{q}{2}}(\tau(\sigma)) M_{2}(T),
\end{equation}
\begin{equation}\label{est:htermW2}
\|\langle z\rangle^{-q}e^{\frac{\alpha
z^{2}}{4}}D_{2}(\sigma)\|_{\infty}\lesssim
\beta^{\frac{1}{4}+\frac{q}{2}}(\tau(\sigma))M_{q}(T),
\end{equation}
\begin{equation}\label{eq:hestSource2}
\|\langle z\rangle^{-q}e^{\frac{\alpha
z^{2}}{4}}D_{3}(\sigma)\|_{\infty} \lesssim
\beta^{1+\frac{q}{2}}(\tau(\sigma))[1+M^{2-p}_{1}(T)+A(T)M_{1}(T)],
\end{equation}
\begin{equation}\label{eq:hnonlinearity2}
\|\langle z\rangle^{-q}e^{\frac{\alpha
z^{2}}{4}}D_{4}(\sigma)\|_{\infty}\lesssim
\beta^{\frac{q}{2}}(\tau(\sigma))M_{q}(T)[M_{q}^{1-p}(T)+M_{q}(T)].
\end{equation}
\end{lemma}
\begin{proof}
The proofs of (~\ref{est:htermW2})-(~\ref{eq:hnonlinearity2}) are
easier than that of Lemma ~\ref{EstDs}, thus is omitted. For
(~\ref{eq:hestF3}) by the fact $q\leq 1$ and the arguments between
(~\ref{eq:compare2}) and (~\ref{eq:compareXiEtaN}) we have
$$
\begin{array}{lll}
\|\langle z\rangle^{-q}e^{\frac{\alpha
z^{2}}{4}}V\eta(\sigma)\|_{\infty}&\lesssim &\|\langle
y\rangle^{-q}\frac{1}{1+\beta(\tau(\sigma))y^{2}}e^{\frac{a(\tau(\sigma))y^{2}}{4}}\xi(y,\tau(\sigma))\|_{\infty}\\
&\leq & \beta^{-\frac{2-q}{2}}(\tau(\sigma))\|\langle
y\rangle^{-2}e^{\frac{a(\tau(\sigma))y^{2}}{4}}\xi(y,\tau(\sigma))\|_{\infty}\\
&\leq &\beta^{\frac{q}{2}}(\tau(\sigma))M_{2}(T)
\end{array}
$$ which is (~\ref{eq:hestF3}).

The proof is complete.
\end{proof}
Rewrite (~\ref{eq:eta4}) to have
$$
P_{1}^{\alpha}\eta(S)=e^{-L_{\alpha}S}P_{1}^{\alpha}\eta(0)+\int_{0}^{S}e^{-L_{\alpha}(S-\sigma)}P_{1}^{\alpha}[-V\eta(\sigma)+\sum_{n=2}^{4}D_{n}]d\sigma,
$$ where, recall the definition
of $S$ in (~\ref{T2}). The propagator estimate of
$e^{-L_{\alpha}\sigma}P_{1}^{\alpha}$ in (~\ref{eq:hPropagator})
yields
\begin{equation}\label{hK123s}
\begin{array}{lll}
\|\langle z\rangle^{-q}e^{\frac{\alpha
z^{2}}{4}}P_{1}^{\alpha}\eta(S)\|_{\infty}\lesssim J_{0}+J_{1}
\end{array}
\end{equation}
where $J_{n}$'s are given by
$$J_{0}:=\beta^{-\frac{q}{2}}(T)e^{[\frac{2}{1-p}-q]\alpha S}\|\langle
z\rangle^{-q}e^{\frac{\alpha z^{2}}{4}}\eta(0)\|_{\infty},$$
$$J_{1}:=\beta^{-\frac{q}{2}}(T)\int_{0}^{S}e^{[\frac{2}{1-p}-q]\alpha(S-\sigma)}[\|\langle z\rangle^{-q}e^{\frac{\alpha z^{2}}{4}}V\eta(\sigma)\|_{\infty}+\sum_{n=2}^{4}\|\langle z\rangle^{-q}e^{\frac{\alpha
z^{2}}{4}}D_{n}(\sigma)\|_{\infty}]d\sigma.$$ We observe that
$\frac{2}{1-p}-q< 0$ by the definition of $q$ and the fact that
$p<-1.$

By the estimates (~\ref{eq:hestF3})-(~\ref{eq:hnonlinearity2}) and
similar procedures as in Sections ~\ref{SEC:EstM1} and
~\ref{SEC:EstM2} we have
$$J_{0}\lesssim M_{q}(0)$$
$$J_{1}\lesssim M_{2}(T)+
\beta^{1/4}(0)[1+M_{q}(T)+M_{1}^{2-p}(T)+M_{1}(T)A(T)]
+M_{q}^{2-p}(T)+M_{q}^{2}(T).$$ These together with the fact that
$\|\langle z\rangle^{-q}e^{\frac{\alpha
z^{2}}{4}}P_{2}^{\alpha}\eta(S)\|_{\infty}=\|\langle
y\rangle^{-q}e^{\frac{a(T) y^{2}}{4}}\xi(T)\|_{\infty}$ and the
definition of $M_{q}$ yield (~\ref{eq:M3}) when $\tau=T$. Since $T$
is arbitrary we have (~\ref{eq:M3}).
%\bibliographystyle{abbrv}
%\addcontentsline{toc}{chapter}{Bibliography}
%\bibliography{biblio}
%\end{document}

\end{document}